\documentclass[11pt,twoside,english]{article}
\usepackage{srcltx,babel,float,amsfonts,amssymb,amsmath,epsfig,amscd}
\usepackage{graphicx,psfrag}
\usepackage[usenames]{color}
\DeclareMathOperator{\im}{Im}
\DeclareMathOperator{\disp}{disp}
\DeclareMathOperator{\dist}{dist}
\DeclareMathOperator{\Diff}{Diff}
\DeclareMathOperator{\diam}{diam}
\DeclareMathOperator{\Id}{Id}
\newtheorem{lemma}{Lemma}[section]
\newtheorem{theorem}{Theorem}
\newtheorem{proposition}[lemma]{Proposition}

\newtheorem{corollary}[theorem]{Corollary}

\newtheorem{definition}[lemma]{Definition}

\newtheorem{remark}[lemma]{Remark}

\newtheorem{assertions}[lemma]{Assertion}

\newtheorem{nada}[lemma]{}
\newtheorem{assertion}[lemma]{Assertion}
\begin{document}
\title{Pesin Entropy Formula for $C^1$ Diffeomorphisms with Dominated Splitting}
\author{Eleonora Catsigeras, Marcelo Cerminara and Heber Enrich
 \thanks{All authors: Instituto de Matem\'{a}tica y Estad\'{\i}stica Rafael Laguardia (IMERL),
 Facultad de Ingenier\'{\i}a,  Universidad de la Rep\'{u}blica,  Uruguay.
 E-mails: eleonora@fing.edu.uy, cerminar@fing.edu.uy,
 enrich@fing.edu.uy. Address: Julio Herrera
  y Reissig 565. Montevideo. Uruguay.}}
\date{September 10th, 2013}

\maketitle

\begin{abstract}
For any $C^1$ diffeomorphism  with dominated splitting we consider a nonempty set of invariant measures which  describes the asymptotic statistics of Lebesgue-almost  all orbits. They are the limits of convergent subsequences of averages of the Dirac delta measures supported on those orbits. We prove that the metric entropy of each of these measures is  bounded from below by the sum of the Lyapunov exponents on the dominating subbundle. As a consequence, if those exponents are non negative, and if the exponents on the dominated subbundle are non positive,      those measures satisfy the Pesin Entropy Formula.
\end{abstract}

\section{Introduction}\label{introduccion}

As pointed out by \cite{pugh,BCS2013} and other authors, there is a gap between the $C^{1+\theta}$ and the $C^1$ Pesin  Theory.  To find  new results that hold for   $C^1$ maps  relatively recent research started     assuming   some uniformly dominated conditions  (see \cite{abc, BCS2013,suntiana, suntian, tahzibi}).

Let us consider    $f \in \Diff ^1(M)$, where $M$  is a compact and connected Riemannian manifold of finite dimension.  We denote by ${\mathcal P} $   the set of all Borel probability measures endowed with the weak$^*$ topology, and   by ${\mathcal P}_f \subset {\mathcal P}$  the  set of $f$-invariant probabilities. We denote by $m$ a normalized  Lebesgue   measure, i.e. $m \in {\mathcal P}$.
For any $ \mu \in {\mathcal P}_f$,  the orbit of   $x$ is   regular for $\mu$-a.e. $x \in M$ (see for instance \cite[Theorem 5.4.1]{barreirapesin}). We denote the Lyapunov exponents of the orbit of $x$  by
$$\chi_1(x) \geq \chi_2 (x) \geq \ldots \geq \chi_{\dim M} (x).  $$
  Let $$\chi^+_i(x) := \max\{\chi_i(x), 0\}. $$

 \noindent \textsc{ Theorem} \textsc{ (Ruelle's Inequality)} \cite{ruelle}  \em

For all $f \in \Diff^1(M)$ and for all $\mu \in {\mathcal P}_f$ \em
 $$h_{\mu} \leq \int \sum_{i= 1}^{\dim M} \chi_i^+  \, d \mu,$$
where $h_{\mu}$ denotes the metric theoretical entropy of $\mu$.

\begin{definition}
Let $f \in \Diff^{1}(M)$ and $\mu \in {\mathcal P}_f$ . We say that
  $\mu$ satisfies the \em Pesin  Entropy Formula\em, and write $\mu \in PF$, if
 $$h_{\mu} = \int \sum_{i= 1}^{\dim M} \chi_i^+  \, d \mu.$$
\end{definition}

  We denote by   $m^u$ the Lebesgue measure along the unstable manifolds of the regular points for which positive Lyapunov exponents and unstable manifolds exist. We denote the (zero dimensional) unstable manifold of $x$ by $\{x\}$, and   in this case we have $m^u= \delta_x$.
For any invariant measure $\mu$ for which local unstable manifolds exist $\mu$-a.e. we denote by $\mu^u$   the conditional measures of $\mu$ along the unstable manifolds, after applying   the local Rohlin decomposition  \cite{ro}.

The following are well known results of the Pesin  Theory under the hypothesis
  $f \in \mbox{Diff }^{2}(M)$:

\vspace{.1cm}

\noindent {\textsc{Pesin  Theorem}} \cite{pesin, ma, barreirapesin}  \em Let $\mu \in {\mathcal P}_ f$ be hyperbolic (namely, $\chi_i(x) \neq 0 $ for all $i$ and for $\mu-$a.e. $x \in M$). If  $\mu \ll m $ then $ \mu^u \ll m^u$  and   $\mu  \in PF$. \em

\vspace{.1cm}

\noindent {\textsc {Ledrappier-Strelcyn-Young Theorem}} \cite{ledrappierstrelcyn, ly}  \em  $\mu \in PF  $ if and only if $ \mu^u \ll m^u$. \em

\vspace{.1cm}

Still in the $C^2$-scenario,  non uniformly partially hyperbolic diffeomorphisms  possess invariant measures $\mu$ such that $\mu^u \ll m^u$; and hence $\mu \in PF$ (see for instance \cite[Theorem 11.16]{BonattiDiazViana}).

\vspace{.2cm}

%This last result states that $\mu$ describes the asymptotic  statistics of a Lebesgue-positive set of orbits (see Remark \ref{asymptoticstatistics}).

\vspace{.3cm}

 The general purpose of this paper is to look for adequate reformulations of some of the above results   which hold  for all
$f \in \mbox{Diff}^{1}(M)$. That is, we would like to know when an invariant measure under $f \in \mbox{Diff}^1(M)$ satisfies Pesin Entropy Formula.

We first recall some definitions and previous results taken from \cite{nuestro}.
 \begin{definition} \em \label{definitionEmpirical} \textsc{ (Asymptotic statistics)}

 Fix $x \in M$. The  \em {sequence  of empirical probabilities} \em of $x$ is $\{\sigma_{n,x}\}_{n \geq 1} \subset{\mathcal P}$, where
      $$ {\displaystyle{ \ \sigma_{n,x} := \frac{1}{n} \sum_{j= 0}^{n-1} \delta_{f^j(x)}}}. $$
      The $p\omega$-\em limit \em of $x$ is
        $$p \omega (x) := \{ \mu \in {\mathcal P}: \  \exists \ n_j \rightarrow + \infty \mbox{ such that } \lim_{j \to \infty} \sigma_{n_j, x} = \mu    \} \subset  {\mathcal P}_f $$
        We say that $p \omega(x)$   describes \em the asymptotic statistics \em of the orbit of $x$.
\end{definition}

\begin{definition} \em \label{definitionBasin} \textsc{ (Basins of statistical  attraction)}

For  any given $\mu \in {\mathcal P} $ the
  {\em {basin of (strong) statistical attraction of $\mu$}}   is
  $$B(\mu) := \left\{ x \in M: \  p \omega(x) = \{\mu\} \right\}. $$
Consider a  metric in ${\mathcal P}$ that induces the weak$^*$ topology and denote it by $\mbox{dist}^*$.

The
\ {\em {basin of { $\mathbf \varepsilon$}-weak statistical attraction of $\mu$}} is
 {  $$B_{\varepsilon}(\mu) := \left\{ x \in M: \  \dist   ^* (p \omega(x), \mu) < \varepsilon   \right\}. $$}
\end{definition}

\begin{definition} \em \label{definitionSRBlike} \textsc{ (SRB, physical and SRB-like measures)}

An invariant probability measure $\mu$ is called {\em SRB } (and we denote $\mu \in $SRB) if  the local unstable  manifolds exist $\mu$- a.e.  and $\mu^u \ll m^u$.

The probability measure $\mu$ is  called  {\em physical } if
$m(B(\mu)) >0.$

If $f \in C^{1+\theta}$  then any hyperbolic   ergodic SRB measure is physical.  Nevertheless,   if $f \in C^1 $, the definition of SRB measure may not be meaningful since there may not exist local unstable manifolds (\cite{pugh, BCS2013}). However, it still makes sense to define when a measure is   physical.

   In the $C^1$-scenario,  we call a  probability measure $\mu$ {\em SRB-like or pseudo-physical} (and we denote $\mu \in $ SRB-like) if
$m(B_{\varepsilon}(\mu)) >0 $ for all $ \varepsilon >0.$

\end{definition}

It is standard to check that the set of SRB-like measures is independent of the metric $\dist^*$ chosen in ${\mathcal P}$ and that it is contained in ${\mathcal P}_f$.

\begin{remark} \label{asymptoticstatistics} \em \textsc{ (Minimal description of the asymptotic statistics of the system)}

Given $f\colon M \to M$, we say that a weak$^*$-compact set ${\mathcal K} \subset {\mathcal P}$ {\em describes the asymptotic statistics} of Lebesgue-almost all orbits of $f$ if $p \omega(x) \subset {\mathcal K}$ for Lebesgue-almost all $x \in M$.
\end{remark}

Theorems 1.3 and 1.5 of \cite{nuestro}  prove that, for any continuous map $f\colon M \mapsto M$ the set of SRB-like measures is nonempty, it contains   $p \omega (x)   $   for Lebesgue-almost all $x \in M$, and it  is the minimal weak$^*$-compact set ${\mathcal K} \subset {\mathcal P}$ such that $p \omega(x) \subset {\mathcal K}$ for Lebesgue-almost all $x \in M$. Therefore,  the set of SRB-like measures \em minimally describes the asymptotic statistics of Lebesgue-almost all   orbits. \em

\vspace{.3cm}

Our focus   is to find     relations, for $C^{1}$ diffeomorphisms, between:

 \noindent $\bullet$  Physical measures and, more generally,  SRB-like  measures.

     \noindent $\bullet$ Invariant measures $\mu$ such that $\mu \in PF$.

Several interesting   results were  already obtained for $f \in \Diff^1 (M)$. First, in \cite{tahzibi} Tahzibi proved the Pesin  Entropy Formula for $C^1$-generic area preserving diffeomorphisms on surfaces. More recently,
Qiu   \cite{qiu} proved that  if   $f$ is a transitive Anosov, then  $C^1$-generically there exists a unique $ \mu  $ satisfying Pesin Entropy Formula.  Moreover $\mu$ is physical and mutually singular with respect to Lebesgue  (cf. \cite{avilabochi}). Finally, we cite:

\noindent     {\textsc {Sun-Tian Theorem}} \cite{suntian}: \em If $f \in \Diff ^1(M) $ has an invariant measure $\mu \ll m$, and if there exists a dominated splitting  $E \oplus F$  $\mu$-a.e. such that  $\chi_{\dim (F)} \geq 0 \geq \chi _{\dim (F) +1} $, then  $    \mu \in PF$. \em

 To prove this theorem Sun and Tian use an approach introduced by Ma\~{n}\'{e} \cite{ma}. In that approach  he gave a new proof of Pesin  Entropy Formula for $f \in \mbox{Diff }^{1 + \theta}(M)$ and hyperbolic $\mu \ll m$. Ma\~{n}\'{e}'s proof     does not directly require the absolute continuity of the invariant foliations. So, it is reasonable to expect that it is adaptable to the $C^1$-scenario.

 We reformulate the technique of Ma\~{n}\'{e} \cite{ma}  to obtain an exact lower bound of the entropy  for  non necessarily conservative $f \in \Diff^1(M)$,  provided that there exists a dominated splitting.
%\vspace{.3cm}

\begin{definition} \textsc{ (Dominated splitting)} \label{definitionSplittingDominado} \em

Let $f\colon M \to M$ be a $C^1$ diffeomorphism on a compact Riemannian manifold.
Let $T M = E \oplus F$ be a continuous and $df$-invariant splitting such that $\mbox{dim}(E), \mbox{dim}(F) \neq 0$. We call $TM = E \oplus F$   a \em dominated splitting \em if there exist  $C >0 $ and $0< \lambda <1$ such that
$${\left\|df^n|_{E_x}\right\|}{\left\|df^{-n}|_{F_{f^n(x)}}\right\|} \leq C \lambda^n, \; \forall x \in M \mbox{ and } n \geq 1.$$

\end{definition}

We will prove   the following results:

\begin{theorem} \label{Teorema2}
   Let $f \in \mbox{Diff}^{\ 1}(M)$ with a dominated splitting $TM= E \oplus F$. Let
    $\mu   $ be an SRB-like measure for $f$.
    Then: \begin{equation}\label{eq0} h_{\mu}(f) \geq \int \sum_{i= 1}^{\dim F} \chi_i \, d \mu.\end{equation}

    \end{theorem}

\begin{corollary} \label{Teorema1}
Under the hypothesis of Theorem \ref{Teorema2}, if  $ \chi_{\dim F} \geq 0 \geq \chi_{\dim F + 1} $, then   $\mu $  satisfies the Pesin  Entropy Formula.

\end{corollary}

 %\begin{theorem}
%  Let $f \in \mbox{Diff}^{\ 1}(M)$ have a dominated  splitting  $TM= E \oplus F$.
%
%   Then:
%
%  \noindent {\em \bf (A)} There exist          measures  $\mu$ characterized by the following properties:
%
%   {\em \bf (A1) } The set of  such measures minimally describe the asymptotic statistics  of Lebesgue-almost all the orbits.
%
%   %\pause
%
%
%    {\em \bf (A2) } They satisfy the inequality:
%
%    \noindent{\em \bf (B)}     Moreover, Property {\em \bf (A1) } implies Inequality {\em \bf  (A2) }.
%
%   \noindent {\em \bf (C)  } If besides $ \chi_{\dim F} \geq 0 \geq \chi_{\dim F + 1} $, then  {any such measure $\mu $  satisfies the Pesin's Entropy Formula.}
%
%
%
%
% \end{theorem}

The proof of Corollary \ref{Teorema1}  is immediate: inequality (\ref{eq0}) and Ruelle's Inequality  imply that $\mu $ satisfies Pesin Entropy Formula. Moreover,   as said in Remark \ref{asymptoticstatistics}, the set of SRB-like measures is nonempty. So, under the hypothesis of Corollary \ref{Teorema1}, there are invariant measures that satisfy the Pesin  Entropy Formula. Besides, they minimally describe  the asymptotic statistics of Lebesgue-almost all  orbits.

Note that according to  Avila and Bochi result   \cite{avilabochi}      the   measures of Theorem \ref{Teorema2} and Corollary \ref{Teorema1} are $C^1$-generically mutually singular with respect to Lebesgue.

\begin{remark}\em

The same arguments of the proof of Theorem \ref{Teorema2}  also work under  hypothesis that are more general than the global dominated splitting assumption. In fact,   if $\Lambda \subset M$ is an invariant and compact topological attractor, and if $\overline V \supset \Lambda$ is a compact neighborhood with dominated splitting $T_{\overline V} = E  \oplus F$, then the same statements and proofs of Theorem  \ref{Teorema2} and Corollary \ref{Teorema1}  hold for $f|_V$.

  Now, let us pose  an example for which   Theorem \ref{Teorema2} and Corollary \ref{Teorema1} do not hold.  Consider the simple eight-figure diffeomorphism in
   \cite[Figure 10.1]{barreirapesin}. In this example, the Dirac-delta measure $\mu$ supported on a fixed hyperbolic point $p$ is physical. Thus $\mu$ is SRB-like. Besides, there exists a dominated splitting $\mu$-a.e. because $p$ is hyperbolic. Nevertheless,  inequality (\ref{eq0}) does not hold because $h_{\mu} = 0$ and the Lyapunov exponent along the unstable subspace of $T_p(M)$ is strictly positive.  So,   the presence of a dominated splitting     \em just    $\mu$-a.e. \em  is  not enough to obtain Theorem \ref{Teorema2}.

\end{remark}

The following question arises from  the statements  of  our results: Does the SRB-like property characterize  all the measures that satisfy Pesin Entropy Formula? The answer is negative. In fact,  the converse statement of Corollary \ref{Teorema1} is false. As a counter-example consider a $C^2$ non transitive uniformly hyperbolic attractor,  with
  a finite set     $ {\mathcal K}= \{\mu_1, \mu_2, \ldots, \mu_k\}$ ($k \geq 2$) of distinct SRB ergodic measures (hence each $\mu_i$ is physical) such that ${\mathcal K}$  statistically attracts Lebesgue-almost every orbit. Therefore,   the set of all SRB-like measures coincides with ${\mathcal K}$ (see Remark \ref{asymptoticstatistics}). So,  $(\mu_1 + \mu_2)/2 \not \in {\mathcal K}$ is not an SRB-like measure. After Corollary \ref{Teorema1},   $\mu_1$ and $ \mu_2 $ satisfy   Pesin Entropy Formula. It is well known that any convex combination of measures that satisfy Pesin Entropy Formula   also satisfies it  (see Theorem 5.3.1 and Lemma 5.2.2. of \cite{ke}). We conclude that $(\mu_1 + \mu_2)/2$ satisfies Pesin Formula but it is not SRB-like.

    \vspace{.3cm}

The   paper is organized as follows:   In Section \ref{sectionReduction} we reduce the proof of Theorem \ref{Teorema2} to    Lemmas \ref{Lemma 1} and \ref{Lemma 2}. In Sections \ref{sectionProofLemma1} and \ref{sectionProofLemma2} we prove Lemmas \ref{Lemma 1} and \ref{Lemma 2} respectively. Finally, in Section \ref{sectionAppendix}   we check some  technical assertions that are used in the proofs of the previous sections.

%    \section{Revisiting the  SRB and SRB-like measures} \label{sectionDefinitions}

%\pause

\section{Reduction of the proof of Theorem \ref{Teorema2}}    \label{sectionReduction}

For the  diffeomorphism $f\colon M \to M$ with dominated splitting $E  \oplus F  = T M$, we denote:
\begin{equation}
\label{eqnPsi}
\psi (x) := - \log \big|\det df(x)|_{F_x}\big| \end{equation}\begin{equation} \label{eqnPsiSubn} \psi_n (x):=   - \log \big|\det df^n(x)|_{F_x}\big| = \sum_{j= 0}^{n-1} \psi \circ f^j(x) = - \log \big|\det df^{-n}(f^n(x))|_{F_{f^n(x)}}\big|
\end{equation}

 Consider a metric $\dist^*$ in the space ${\mathcal P}$ of all Borel probability measures  inducing its weak$^*$ topology. For all $\mu \in {\mathcal P}$, for all $\varepsilon >0$ and for all $n \geq 1$, we denote:
\begin{equation} \label{equationCn}
C_n(\varepsilon) := \{x \in M: \ \dist^* (\sigma_{n,x}, \mu) < \varepsilon\},\end{equation}
where $\sigma_{n,x}$ is the empirical probability   according to Definition \ref{definitionEmpirical}. We call $C_n(\varepsilon)$ the \em  approximation  up to time $n$ \em   of the  basin  $B_{\varepsilon}(\mu)$ of $\varepsilon$-weak statistical attraction   of the  measure $\mu$ (cf. Definition \ref{definitionBasin}).
\begin{proposition} \label{theoremNuevo}
Let \em $f \in \mbox{Diff} ^1(M)$ \em with a dominated splitting $TM= E \oplus F$. There exists a weak$^*$ metric $\dist^*$ in ${\mathcal P}$, such that for any $f$-invariant probability measure $\mu$ the following inequality holds:
\begin{equation} \label{eqnNueva}\lim_{\varepsilon  \rightarrow 0^+} \limsup_{n \rightarrow + \infty} \frac{\log m(C_n(\varepsilon))}{n} \leq h_{\mu}(f) + \int \psi \, d \mu,\end{equation}
where $m$ is the Lebesgue measure.
\end{proposition}

\vspace{.5cm}

   We note that the term  $h_{\mu}(f) + \int \psi \, d \mu$ is non negative due to Ruelle's Inequality.   Nevertheless, it is   bounded from below by  inequality (\ref{eqnNueva}), which relates it with the Lebesgue measure $m$.

\vspace{.3cm}

\noindent At the end of this section, we reduce the proof of Proposition \ref{theoremNuevo}  to Lemmas \ref{Lemma 1} and \ref{Lemma 2}.  Along the remaining sections   we prove these two lemmas. Now, let us prove the following assertion:

\vspace{.2cm}

\noindent \textsc{ Proposition \ref{theoremNuevo} implies Theorem \ref{Teorema2}.}

\vspace{.2cm}

{\em Proof: }

 Let $\mu$ be $f$-invariant. Assume that $\mu$ does not satisfy inequality (\ref{eq0}). In other words,
 $$ h_{\mu}(f) +   \int \psi \, d \mu = - r < 0.  $$
 From Proposition \ref{theoremNuevo}, for all $\varepsilon >0$ small enough there exists $N \geq 1$ such that
 $$\frac{\log m(C_{n}(\varepsilon))}{n} \leq \frac{-r}{2} \ \ \forall \ n \geq N.$$

 Since $r >0$, we deduce that
 $\sum_{n= 1}^{+ \infty} m(C_{n}(\varepsilon)) < + \infty.$
 Thus, by Borel-Cantelli Lemma   the set
$\bigcap_{N \geq 1} \bigcup_{n \geq N} C_{n}(\varepsilon)$
has zero $m$-measure. By Definition \ref{definitionBasin} we have
$B_{\varepsilon}(\mu) \subset \bigcap_{N \geq 1} \bigcup_{n \geq N} C_{n}(\varepsilon).$
So,   $m(B_{\varepsilon}(\mu))= 0$, and applying Definition \ref{definitionSRBlike} we conclude that $\mu$ is not SRB-like, proving Theorem \ref{Teorema2}.
 \hfill $\Box$

\vspace{.2cm}

\noindent \textsc{ Proposition \ref{theoremNuevo} follows from Lemmas \ref{Lemma 1} and \ref{Lemma 2}.}

\vspace{.2cm}

 To prove Proposition \ref{theoremNuevo},  we   take from  \cite{suntian} the idea of using Ma\~{n}\'{e}'s approach \cite{ma}. Nevertheless, we use this approach in a   distinct context (i.e. we do not assume $\mu \ll m$) and   apply different arguments.
In \cite{ma} Ma\~{n}\'{e} considers $f \in C^{1+\theta}$ and constructs a $C^1$ foliation ${\mathcal L}$, which is not necessarily invariant, but approximates the unstable invariant foliation.
On the one hand, the given invariant measure $\mu \ll m$ has absolutely continuous conditional measures along the leaves of ${\mathcal L}$, because $m$ has. On the other hand, the hypothesis $f \in C^{1+\theta}$ allows Ma\~{n}\'{e} to use   the Bounded Distortion Lemma. So, he  obtained Pesin Entropy Formula after taking  $f^n {\mathcal L}  $ convergent to the unstable foliation.

In our case these arguments fail to work,  except one. There still exists a $C^1$  (non invariant) foliation ${\mathcal L}$ whose tangent sub-bundle approximates the dominating sub-bundle $F$. Besides, since ${\mathcal L}$ is $C^1$, the conditional   measures of $m$ (not of $\mu$) along the leaves of   ${\mathcal L}$ are absolutely continuous. But we   have neither the hypothesis   $\mu \ll m$ nor the $C^{1+\theta}$ regularity of $f$. Also   an invariant foliation  to which $f^n {\mathcal L}$ would converge, may fail to exist.  The role of the following   Lemmas \ref{Lemma 1} and \ref{Lemma 2}  is to overcome   these problems. Before stating them,  we adopt the following:

 \noindent\textsc{ Notation. } Let ${\mathcal B}$ be the Borel $\sigma$-algebra on  the manifold $M$. We denote by $\alpha = \{X_i\}_{1 \leq i \leq k} $ a finite partition of $M$, namely:
 \begin{center}

   \noindent $X_i \in {\mathcal B} $ for all $ \ 1 \leq i \leq k$,

   \noindent  $X_i\bigcap X_j = \emptyset$ if $i \neq j$,

   \noindent  $\bigcup_{i=1}^k X_i = M$.

 \end{center}
   We write $f^{-j}(\alpha) = \{f^{-j}(X_i)\}_{1 \leq i \leq k}$.

For any pair of finite partitions $\alpha = \{X_i\}_{1 \leq i \leq k}$ and $\beta = \{Y_j\}_{1 \leq j \leq h}$ we denote $$\alpha \vee \beta = \{X_i \cap Y_j: 1 \leq i \leq k,\ 1 \leq j \leq h, \  X_i \cap Y_j \neq \emptyset\},$$ $$\alpha^n = \bigvee_{j= 0}^{n} f^{-j}(\alpha).$$

 \begin{lemma}
  \label{Lemma 1}
  \textsc{ {(Upper bound of the Lebesgue measure $ {m}$)}}

   For all   $ \varepsilon >0 $ there exists $   \delta >0$ such that for every
   finite partition $\alpha $ with $\diam (\alpha) < \delta$
there exist a sequence $\{\nu_n\}_{n \geq 0}$ of finite measures and a constant $K >0$ such that:

    {\em \bf (i)} $\nu_n(X) < K $ for all $ \ X \in \alpha^n = \bigvee_{j= 0}^{n} f^{-j}(\alpha) $, for all $   \ n \geq 0$.

    {\em \bf (ii)} The following inequality holds for all $  C \in {\mathcal{B}}  $ and for all $ n \in \mathbb{N}^+:$
     $$  \ m(C) \leq K e^{n \varepsilon} I(\psi_n, C, \nu_n), \ \ \mbox{ where }  $$
    \begin{equation}
     \label{eqn*}
     I(\psi_n, C, \nu_n) : = \int _{C} e^{\psi_n} \, d \nu_n.\end{equation}
\end{lemma}

We will prove Lemma \ref{Lemma 1} in Section \ref{sectionProofLemma1}.

Before stating the second lemma, recall equality (\ref{equationCn}).

\begin{lemma}
\label{Lemma 2}
\textsc{ (Lower bound of the metric entropy)}

   There exists a metric $\dist^*$ in ${\mathcal P}$ with the following property:

  For all  $  \mu \in {\mathcal P}_ f$ and for all $   \varepsilon, \ \delta >0 $ there exist a finite partition $\alpha$ satisfying $\diam  \alpha < \delta$,  and a real number $\varepsilon_0^* >0$ such that:

   For all $0 <\varepsilon^*< \varepsilon^*_0$, and
for any sequence $\{\nu_n\}_{n \geq 0}$ of finite measures such that   there exists $   K >0$ satisfying $\nu_n(X) < K $ for all $ \ X \in {\alpha}^n $ for all $ n \geq 0$, the following inequality holds:
   $$\limsup_{n \rightarrow + \infty} \frac{1}{n} \log I(\psi_n, C_n(\varepsilon^*), \nu_n) \leq \varepsilon + h_{\mu} ( \alpha) + \int \psi \, d \mu.$$

\end{lemma}

We will prove Lemma \ref{Lemma 2} in Section \ref{sectionProofLemma2}.

\vspace{.3cm}

  To end this section let us prove that
Lemmas \ref{Lemma 1} and \ref{Lemma 2} imply  Proposition \ref{theoremNuevo}:

\vspace{.2cm}

 {\em Proof: }
  { Let $\mu \in {\mathcal P}_f$  and $\varepsilon >0$.
Consider $\delta >0$ obtained from   Lemma \ref{Lemma 1}. Applying Lemma
     \ref{Lemma 2}, construct the partition $\alpha$, the number $\varepsilon^*_0$ and the sequence   $\{C_n(\varepsilon^*)\}_{n \geq 0} \subset {\mathcal B}$ for any $0 <\varepsilon^*< \varepsilon^*_0$.

     Apply again   Lemma \ref{Lemma 1}  to obtain the sequence $\{\nu_n\}_{n \geq 0}$  of finite measures and the constant $K >0$.

We now apply again Lemma \ref{Lemma 2} to deduce: \begin{equation}\limsup_{n\to \infty } \frac{1}{n} \log I(\psi_n, C_n(\varepsilon^*), \nu_n) \leq \varepsilon + h_{\mu} ( \alpha) + \int \psi \, d \mu  \ \ \forall \ 0 <\varepsilon^* < \varepsilon^*_0\label{equation1}
     \end{equation}

    Besides, by  Lemma \ref{Lemma 1}:
    \begin{equation}
     \frac{1}{n} \log m(C_n(\varepsilon^*)) \leq \frac{\log K}{n} + \varepsilon +\frac1{n} \log I(\psi_n, C_n(\epsilon^*), \nu_n). \label{equation2}
     \end{equation}

 We join the two inequalities (\ref{equation1}) and (\ref{equation2}) to deduce that:
    $$\limsup_{n \to \infty}  \frac{1}{n} \log m(C_n(\varepsilon^*)) \leq 2 \, \varepsilon + h_{\mu} ( \alpha) + \int \psi \, d \mu  \ \ \forall \ 0 < \varepsilon^* < \varepsilon^*_0.$$

Taking $\varepsilon^* \rightarrow 0^+$ we obtain:
$$\lim_{\varepsilon^* \rightarrow 0^+} \limsup_{n \to \infty}  \frac{1}{n} \log m(C_n(\varepsilon^*)) \leq 2 \, \varepsilon + h_{\mu} ( \alpha) + \int \psi \, d \mu.$$
     Since $\varepsilon>0$ is arbitrary, we deduce inequality (\ref{eqnNueva}), as wanted.
       \hfill $ \Box$}

\section{Proof of Lemma \ref{Lemma 1}} \label{sectionProofLemma1}

To prove Lemma \ref{Lemma 1} we will use  the technique of the dispersion  of  Hadamard graphs, following Ma\~{n}\'{e} in \cite{ma}.

 \noindent \textsc{ Notation:} First take a fixed value of $\delta >0$ small enough  such that  $\exp^{-1}_x$ is a diffeomorphism  from $B_{3 \delta}(x)$ onto its image in $T_xM$ for all $  x \in M$. Fix $x \in M$.
  Denote
  $${\mathbf B}^{E_x} _{\delta} ({\mathbf 0}) := \{v \in E_x: \|v\| \leq \delta\}, $$
  $${\mathbf B}^{F_x} _{\delta} ({\mathbf 0}) := \{v \in F_x: \|v\| \leq \delta\}, $$
  $$
 {\mathbf B}^{T_xM}_{\delta}({\mathbf 0}) := {\mathbf B}^{E_x} _{\delta} ({\mathbf 0}) \oplus {\mathbf B}^{F_x} _{\delta} ({\mathbf 0})    .$$

%\pause
Denote by $\pi_{E_x}$ (resp. $\pi_{F_x}$) the projection of $T_xM$ on $E_x$ along $F_x$ (resp. on $F_x$ along $E_x$), and $\ \ \gamma := \max_{x \in M} \{\|\pi_{E_x}\|, \|\pi_{F_x}\|\}$.

  For any $   v \in {\mathbf B}^{T_xM}_{\delta}({\mathbf 0})$ we denote   $v_1:=\pi_{E_x} v , \ \ v_2:=\pi_{F_x} v.$
%\pause

%\vspace{-.15cm}

  \begin{definition} \label{definitiongraph} \em $G$ is a \em Hadamard graph \em (or simply \lq\lq a graph\rq\rq) if

        $ G\colon {\mathbf B}_{\delta} ^{E_x}({\mathbf 0}) \times {\mathbf B}_{\delta}^{F_x}({\mathbf 0})  \to {\mathbf B}_{\delta}^{E_x}({\mathbf 0})$,

        $ G(v_1, 0) = 0 $ for all $ v_1 \in {\mathbf B}_{\delta} ^{E_x}({\mathbf 0})$  and

  $\Phi (v_1, v_2) = v_1 + v_2 + G(v_1, v_2) \in {\mathbf B}^{T_xM}_{2 \delta}({\mathbf 0})$  is a $C^1$-diffeomorphism onto its image. (See Figure \ref{Figure}.)
  \end{definition}

\begin{center}

%\vspace{-.25cm}
{\begin{figure}
[h]
\begin{center}\includegraphics[scale=.6]{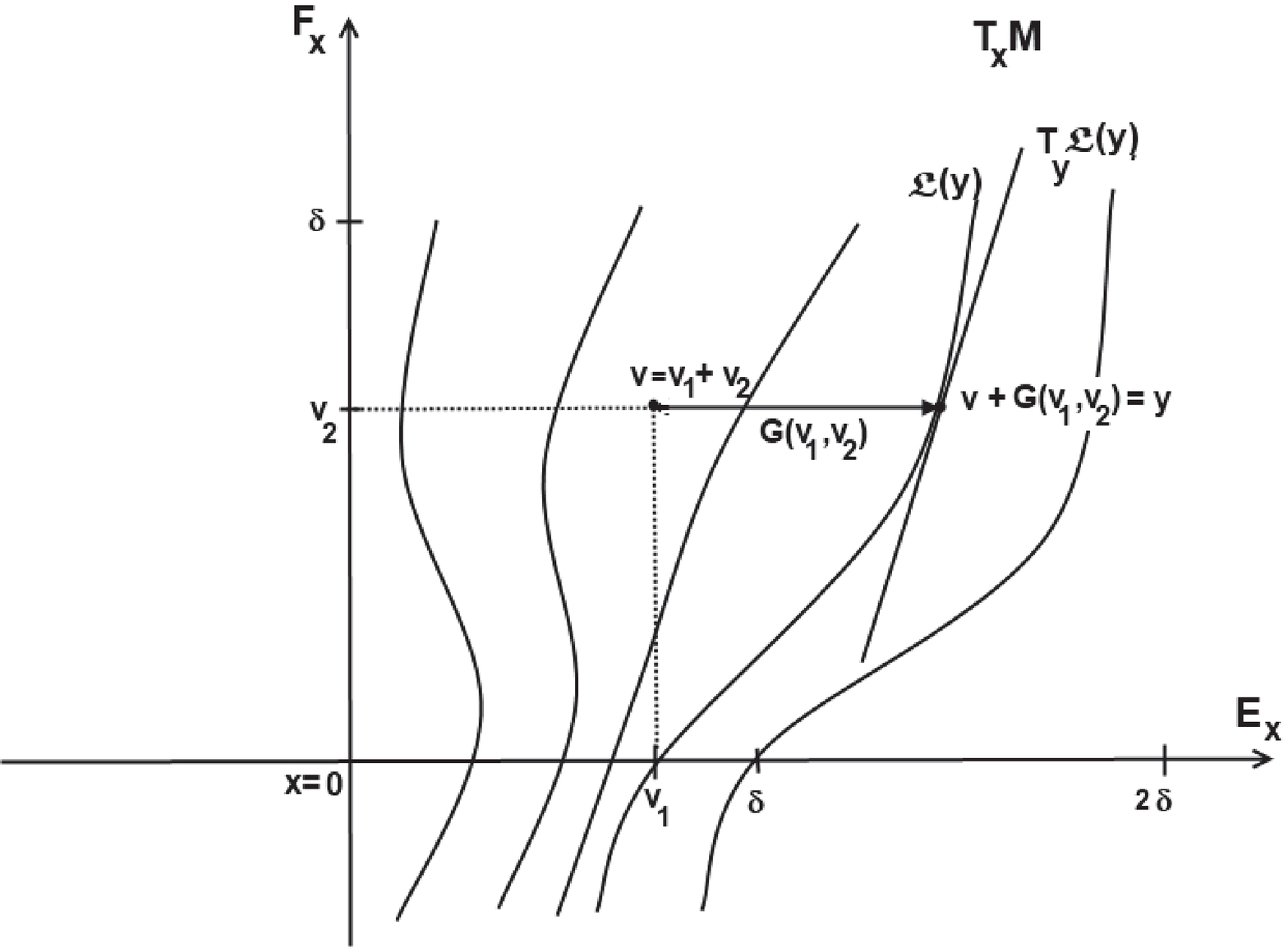}
\caption{\label{Figure} The foliation ${\mathcal L}$ associated to a Hadarmard graph (We omit the exponential map $\exp_x$)}
\end{center}
\end{figure}}
\end{center}

  The  foliation ${\mathcal L}$ associated to the graph $G$ is the foliation whose leaves   are   parametrized on $v_2 \in {\mathbf B}_{\delta}^{F_x}({\mathbf 0}) \subset F_x$, with constant $v_1$, by the diffeomorphism:
  $$\exp_{x} (\Phi(v_1, \cdot)) = \exp_x (v_1 + \cdot + G(v_1, \cdot))  $$
   In Figure \ref{Figure} we draw the foliation ${\mathcal L}$. To simplify the notation we omit the exponential map $\exp_x$ and  denote $y = v_1 + v_2 + G(v_1, v_2)$. The leaf containing $y$ is denoted by ${\mathcal L}(y)$.

  \begin{definition}
  \label{DefinitionDispersion} \em
   \textsc{ Dispersion of $\mathbf{G}$}

   The dispersion of the graph $G$ is

    $$\disp   G  = \max_{v \in {\mathbf B}^{T_xM}_{\delta}({\mathbf 0})} \left\{\left\|\frac{\partial G}{\partial v_2}(v_1, v_2) \right\|\right\},$$
   where $v=v_1+v_2$ and {$\partial G/\partial v_2$ denotes the Fr\'{e}chet derivative of   $$G(v_1, \cdot)\colon {\mathbf B}_{\delta}^{F_x}({\mathbf 0}) \to {\mathbf B}_{\delta}^{E_x}({\mathbf 0})$$ with a constant value of $v_1 \in {\mathbf B}_{\delta}^{E_x}({\mathbf 0})$.}

 \end{definition}

We   denote by $m^{W}$ the Lebesgue measure along an embedded local submanifold $W \subset M$.

 \begin{assertions}
\begin{equation}  \label{equation5} T_y {\mathcal L}(y) = \left(\Id |_{F_x} + \frac{\partial G}{\partial v_2} \right) F_x,\end{equation}
\begin{equation} \label{ocho} m^{{\mathcal L}(y)}({\mathcal L}(y))  \leq [(1 + \disp  \, G) \delta ]^{\dim F}, \end{equation}
For all $  \ \varepsilon >0 $ there exists $ c >0$ such that, if $\disp \, G \leq c $, then $$ \dist   (T_y {\mathcal L}(y), F_x) <  \frac{\varepsilon}2 \ \ \ \ \forall \ y \in \im(\Phi). $$

For  such a value of $c >0$ (depending on $\varepsilon >0$), there exists $ \delta_1 >0 $ such that, if  $ \dist  \,(x,y) < \delta_1 $, then $ \ y \in \im (\Phi) $ and \begin{equation}\label{9} \dist   (T_y{\mathcal L}(y), F_y) < \varepsilon \end{equation}
\end{assertions}

{\em Proof:} The   assertion  follows from   the properties that were established in the definition   of  Hadamard graphs and their associated foliations, and from the definition of dispersion. In particular (\ref{9}) holds because of the continuous dependence of the splitting $E_y \oplus F_y$ on the point $y$. \hfill $\Box$

\newpage

  \begin{nada}
  \textsc{ Iterating  the local foliation ${\mathcal L}$}
  \end{nada}

Denote by $B^n_{\delta}(x)$ the dynamical ball defined as $$B^n_{\delta}(x):= \{y \in M: \dist  (f^j(x), f^j(y)) < \delta \ \forall\, 0 \leq j \leq n\}.$$
Take any graph $G$ in ${\mathbf B}^{T_xM}_{\delta}({\mathbf 0})  $ such that $\disp  G < 1/2$, and consider its associated local foliation ${\mathcal L}$.  Construct the   image  $f^n({\mathcal L})$ in the dynamical ball $B^n_{\delta}(x)$, i.e.:
    $$f^n({\mathcal L} \cap B^n_{\delta}(x))  = f^n \exp_x (v_1   \!+ v_2   +  G(v_1, v_2)) $$
 for all $(v_1, v_2)  \in {\mathbf B}^{E_x}_{\delta}({\mathbf 0})\times {\mathbf B}^{F_x}_{\delta}({\mathbf 0})$  such that $ \exp_x(v_1  +   v_2   +  G(v_1, v_2)) \in B^n_{\delta}(x).$

\begin{lemma} \textsc{ (Reformulation of Lemma 4 of \textsc{ \cite{ma}})}  \label{LemmaMane}

There exists  $ 0 <c'<  1/2$ depending only on $f$, such that for all $0 <c < c' $ there are $\delta_0,\,  n_0 >0$ such that for any point $x \in M$,  if ${\mathcal L}$ is the local foliation associated to a  graph $G$ defined on $T_xM$ with
$$\disp  \, G < c ,$$      then for all $n \geq 0$ the iterated foliation $f^n({\mathcal L} \cap B^n_{\delta_0}(x))$ is contained in the associated foliation of  a  graph $G_{n}$ defined on $T_{f^n (x)}M$, and \em
\begin{equation}  \label{eqnLemmaMane} \disp  \, G_{n } < c \mbox{ for all } n \geq n_0.\end{equation}
\em Besides, for all $y \in B_{\delta_0}^n(x)$ the image $f^n({\mathcal L}(y) \cap B^n_{\delta_0}(x) )$ is contained in a single leaf of the foliation associated to $G_{n}$.

   \end{lemma}

   \noindent {\em Proof: }

  \noindent\textsc{ Step 1.}    Choose $n_0 \geq 0$  such that \begin{equation} \label{eqn22} \big\|df^{n }|_{E_z} \big \| \, \, \big\|df^{-n }|_{F_{f^{n }(z)}}  \big\| < 1 \ \ \forall \  n \geq n_0, \ \ \forall \ x \in M.\end{equation}

  For such a fixed value of $n_0$,  take $\delta_0>0$ so that for all $x \in M$, for all $0 \leq n \leq n_0$, and for any graph  $G$ defined in ${\mathbf B}_{\delta}^{T_xM}({\mathbf 0}) $ with $\disp G  < 1/2$, there exists a graph $   G_{n}$ defined on  $\mathbf{B}_{\delta}^{T_{f^n(x)} M}({\mathbf 0})  $ satisfying the following condition:

\begin{equation}
\label{eqn20}\begin{array}{l}
               \mbox{ for any } \ y = \exp_{x} (v_1 + v_2 + G(v_1, v_2)) \in B_{\delta_0}(x)  \\
               \mbox{ there exists }
               (u_1, u_2)  \in {\mathbf B}_{\delta}^{E_{f^n(x)}}({\mathbf 0})\times {\mathbf B}_{\delta}^{ F_{f^n(x)}}({\mathbf 0}) \\
               \mbox{ where } u_1 \mbox{ depends only on } v_1 \mbox{ and } \\
                 f^n(y)= \exp_{f^n(x)}(  u_1 \! \!+  \!u_2 \! + \!  G_{n} \! (u_1, u_2)).
             \end{array}
\end{equation}
  %\begin{equation}
%   \label{eqn20a}
%   \mbox{ For any } \ y = \exp_{x} (v_1 + v_2 + G(v_1, v_2)) \in B_{\delta_0}(x) \mbox{ there exists } $$ $$ (u_1, u_2)  \in {\mathbf B}_{\delta}^{E_{f^n(x)}}({\mathbf 0})\times {\mathbf B}_{\delta}^{ F_{f^n(x)}}({\mathbf 0}), \end{equation}    \begin{equation}\label{eqn20b} \mbox{ where } u_1 \mbox{ depends only on } v_1 \end{equation} \begin{equation}  \end{equation}

In  Assertion \ref{assertionG*1} of the appendix  we show that such   $\delta_0 >0$ exists. We note that     the above assertion is true for any initial graph $G$ with dispersion smaller than $1/2$ and that   Statement  (\ref{eqn20}) a priori \em only \em holds   if $0 \leq n \leq n_0$. The assertion that $u_1$ depends only on $v_1$ implies that the image   $f^n(z)$ of any point $z \in B_{\delta_0} (x)$  in the   leaf ${\mathcal L}(y)$ associated to the graph $G$, is contained in   the   leaf of $f^n(y)$ associated to the graph $G_n$.

 \noindent\textsc{ Step 2.} With $\delta_0>0$ fixed as above, there exists $0 < c' < 1/2$ such that
for any graph  $ G$ with $\disp  \, G <   c'$ and for all $ n \geq 0$,
if $G_n$ is the graph defined in ${\mathbf B}^{T_{f^n(x)} M}_{\delta} ({\mathbf 0}) $ satisfying   (\ref{eqn20}),   then:
   \begin{equation}
     \label{eqn21}
     \left\| \displaystyle{\frac{\partial  G _{n}}{\partial u_2}(u_1, u_2)} \right\|     \leq    \|df^n|_{E_x} \| \cdot     \disp  G  \cdot    \big\| df^{-n}|_{F_{f^n(x)}} \big\|    \ \ \forall \ y \in B_{\delta_{ 0}}^{n}(x) \end{equation}
We prove this  statement in Assertion \ref{assertionG*2} of the appendix.

\noindent\textsc{ Step 3.} Due to the construction of $\delta_0$ in Step 1, inequality (\ref{eqn21}) holds   in particular  for  $n= n_0$ for any $G$ such that $\disp  G  < c'$. Therefore, using Inequalities (\ref{eqn22}) and (\ref{eqn21})  and Definition \ref{DefinitionDispersion},  we obtain:
    $$ \disp  G_{n_0} \leq \disp G <   c' \ \ \forall \ G \mbox{ s.t. } \disp  G < c',$$
Moreover, if $\disp G < c < c'$, then $\disp G_{n_0} < c < c'.$

\noindent\textsc{ Step 4.} From the construction of $\delta_0$ in   Step 1 and using that $\disp  G_{n_0} < c < c' < 1/2$, we deduce that the graph $G_n$  exists  for all $n_0 \leq n \leq 2n_0$. Moreover, $\|df^n|_{E_x}\| \cdot \|df^{-n}|_{F_{f^n(x)}}\| < 1$ for all $n \geq n_0$. So, applying inequality  (\ref{eqn21})   we obtain $\disp  G _n < c$     for all $n_0 \leq n \leq 2n_0$. Finally,     applying inductively    Assertions (\ref{eqn22}) and (\ref{eqn20})  we conclude  that   the graph $G_n$ exists for all $n \geq 0$ and
$\disp  G_n < c$ for all $n \geq n_0$.
\hfill $\Box$

\vspace{.3cm}
Once  the constant $c'$ of Lemma \ref{LemmaMane} is fixed, depending only on $f$, one obtains the following property that allows to move the reference point $x$ (used to construct the  graph $G$ on ${\mathbf B}^{T_xM}_{\delta}({\mathbf 0})  $), preserving the \em same \em   associated local foliation ${\mathcal L}$ and the        uniformity of the upper bound of its dispersion:
\begin{lemma}
\label{lemmaNuevoSet2012} \em

 \em For all $0 < c < c'$ there exists $\delta_1 >0$ such that, for any $x \in M$ and for any graph $G$ with   $\disp(G) < c/2$ defined in ${\mathbf B}^{T_xM}_{\delta}({\mathbf 0})  $, the associated foliation ${\mathcal L}$   in the neighborhood $B_{\delta_1}(x)$ is also associated to a graph $G'$ defined in ${\mathbf B}^{T_zM}_{\delta}({\mathbf 0})$   for any $ z \in B_{\delta_1}(x)$.
 Besides  $\disp (G') < c$.
\end{lemma}

{\em Proof: }
The splitting $E_z \oplus F_z$ depends continuously on $z \in M$. Then    $\pi_{E_z}$ and $\pi_{F_z}$ also depend continuously on $z$. Therefore, for all $\varepsilon >0 $ there exists $\delta_1>0$ such that $$\|\pi_{F_x}|_{E_z}\| < \varepsilon, \ \ \ \ \|\pi_{E_x}|_{F_z}\| < \varepsilon \ \ \ \mbox{ if } \dist (x,z) < \delta_1. $$
(For simplicity in the notation  in  the above inequalities we omit the derivative of $\exp^{-1}_{x} \circ \exp_z$  which identifies $T_z M$ with $T_xM$.)

We claim  that if $\delta_1 >0$ is small enough   then, for any graph $G$   defined on ${\mathbf B}^{T_xM}_{\delta}({\mathbf 0})  $, and for any point $z$  such that $\dist (z,x) <  \delta_1  $,   there exists a graph $G'$  defined on ${\mathbf B}^{T_zM}_{\delta}({\mathbf 0})  $ such that the local foliations associated to $G$ and $G'$ coincide in an open set where both are defined. In fact,   $G'$ should satisfy the following equations:
\begin{equation} \label{eqn31c} u_1 + u_2 + G'(u_1, u_2) = v_1 + v_2 + G(v_1, v_2), \end{equation}
$$u_1, G'(u_1, u_2) \in E_z, \ \ \ \ \ u_2 \in F_z, \ \ \ \ \ G'(u_1, 0) = 0.$$
Since by hypothesis $G$ is a graph, it is $C^1$ and
$$v_1, \ \ G(v_1, v_2) \in E_x, \ \ v_2 \in F_x,  \ \ \ G(v_1, 0) = 0.$$
The above equations are solved by
\begin{equation} \label{eqn31a}u_1 := \pi_{E_z} (v_1),\ \ \ \ \ u_2 := \pi_{F_z} (v_1 + v_2 + G(v_1, v_2)),\end{equation}
\begin{equation} \label{eqn31b}G' = -u_1 + \pi_{E_z} (v_1 + v_2 + G(v_1, v_2)). \end{equation}
The   two equalities in (\ref{eqn31a})  define a local diffeomorphism $\Psi (v_1, v_2) = (u_1, u_2)$. In fact, on the one hand $u_1 =$ $\pi_{E_z}|_{E_x}(v_1)$, where $\pi_{E_z}|_{E_x}$ is a diffeomorphism (which is linear and uniformly near the identity map, independently of the graph $G$). On the other hand, for $v_1$ constant, the derivative with respect to $v_2$ of $\pi_{F_z}(v_1 + v_2 + G(v_1, v_2)) $ is $\pi_{F_z}|_{T_y{\mathcal L}(y)}$, which is, independently of the graph $G$, uniformly near $\pi_{F_x}|_{T_y{\mathcal L}(y)} = Id|_{F_x}$. Thus, $\Psi$ is a local diffeomorphism    $C^1$ near the identity map provided that $\delta_1$ is chosen small enough (independently of the given graph $G$).

From the above construction  we deduce that  the composition of the  mapping $\Psi (u_1, u_2) = (v_1, v_2)$ with the mapping $(v_1, v_2) \mapsto G'$ defined by (\ref{eqn31b}),  is of $C^1$ class. Therefore $G'(u_1, u_2)$ is $C^1$ dependent on $(u_1, u_2)$. Besides $G'(u_1, 0)= 0$ because $G(v_1, 0) = 0$. Due to Identity (\ref{eqn31c}), the application $\phi'$ defined by $\Phi'(u_1, u_2) := u_1 + u_2 + G'(u_1, u_2)$ coincides with the application $\Phi(v_1, v_2) := v_1 + v_2 + G(v_1, v_2)$.

 Due to Definition \ref{definitiongraph} the mapping $\Phi$  is a local diffeomorphism. So $\Phi'$ is also a local diffeomorphism. Thus $G'$ satisfies  Definition \ref{definitiongraph} of Hadamard graph. The first claim is proved.

The diffeomorphism $\Psi (v_1, v_2) = (u_1, u_2)$ as constructed above, converges   to the identity map in the $C^1$ topology, when $\delta_1 \rightarrow 0 ^+$, and  uniformly for all graphs  $G$ defined in ${\mathbf B}^{T_xM}_{\delta}({\mathbf 0})  $. Thus, by Identity (\ref{eqn31c}), $\|G' - G\|_{C^1}$ converges uniformly to zero, independently of the given graph $G$, when $\delta_1 \rightarrow 0$. This implies, in particular, that $\partial G'(u_1, u_2)/ \partial u_2$ converges uniformly to $\partial G(v_1, v_2) / \partial v_2$ when $\delta _1  \rightarrow 0.$ Thus, for any constant $c/2 >0$ there exists $\delta_1>0$, which is independent of the graph $G$,  such that $| \disp (G') - \disp (G)| < c/2$. In other words, $\disp (G') < c$ for all $G$ with $\disp(G) < c/2$, as   wanted.
 \hfill $\Box$

\vspace{.3cm}

We are   ready to prove the following Proposition, for all $f \in \Diff ^1(M)$ with a dominated splitting $TM= E \oplus F$.

\begin{proposition}  \label{PropositionMane}

 For all { $  \varepsilon > 0 $ there are $    \delta_0,  K,   n_0 >0,  $ and  a finite family of local foliations ${\mathcal L}$, each one defined in an open ball  of a  given finite covering of $M$ with $\delta_0$-balls,   such that:

%\pause

\noindent {\rm (a)}    ${\mathcal L}  $ is $C^1$- trivializable and its leaves are $\dim F$-dimensional,

%\vspace{.1cm}

%\pause

\noindent{\rm (b)}   \em $\dist  \left(F_{f^n(x)}, T_{f^n(x)} f^n ({\mathcal L} (x))\right) < \varepsilon $ for all $ x$ and for all $\ n \geq n_0$, \em

%\vspace{.2cm}
%\pause

%\noindent\textsc{ (c){\rm\&}(d)} \ $\forall \ n \geq 0$ and $\forall \ x,y $ such that $ y  \in B_{\delta_0}^n(x)$ \em (dynamical ball): \em

%\pause
%\vspace{.1cm}

\noindent{\rm (c)}    the following assertion holds for all $ n \geq 0$ and for all $  x,y $ such that $ y  \in B_{\delta_0}^n(x)$: $$m^{f^n({\mathcal L} (y))}(f^n ({\mathcal L} (y) \cap B^n_{\delta_0}(x))) \leq K, $$

%\vspace{.1cm}

%\pause

\noindent{\rm (d)}  the following inequality holds  for all $  n \geq 0$ and for all $  x \in M$:  $${  {e^{-n \varepsilon}}{K}^{-1} \leq \frac{\left|\det df_x^{n}|_{ T_{x}({\mathcal L}(x))}\right|}{\left|\det df_x^{n}|_{F_x} \right|} \leq K e^{n \varepsilon}}.$$

}
\end{proposition}
{\em Proof: }
Consider the constant $c'$ determined by Lemma \ref{LemmaMane}. For each point $x \in M$ construct a local foliation ${\mathcal L}$ from a graph $G$ defined on $T_xM$, with dispersion smaller than a constant $c/2$ such that $0 <c < c' < 1/2$. The  constant $c$ will be fixed later taking into account the given value of $\varepsilon>0$.

After Lemma \ref{lemmaNuevoSet2012}, there exists $\delta_1 >0$ such that, for all $x \in M$  the graph $G$ defined on ${\mathbf B}^{T_xM}_{\delta}({\mathbf 0})  $  is redefined on ${\mathbf B}^{T_zM}_{\delta}({\mathbf 0})  $,  for any point $z \in B_{\delta_1}(x)$, preserving the same associated foliation and having dispersion upper bounded by $c$.  Fix  $\delta_0, n_0$ (depending on $c$) by    Lemma~\ref{LemmaMane} and such that $\delta_0 < \delta_1$.
For any given finite covering of $M$ with balls $B_{\delta_0}(x_i)$, fix a finite family $\{\mathcal L_i\}_{1 \leq i \leq k}$ of local foliations so constructed, one in each ball of the covering.

 By the definition of graph, each foliation ${\mathcal L}$ of the finite family constructed above, is $C^1$-trivializable and its leaves have the same dimension as the dominating subbundle $F$. Thus Assertion {\rm (a)}  is proved.

%\noindent Now, we will prove {\rm (c)}
%From Inequality (\ref{ocho}), and for any value of $c >0$, we deduce that if
%the graph $G$ has dispersion smaller than $c$ then:
%\begin{equation} \label{eqC} {  m^{{\mathcal L}(y)} ({\mathcal L}(y)) \leq [(1+ c) \delta]^{\dim F}     } \end{equation}

From inequality (\ref{9}),  given   $    \varepsilon' >0 $ (a fixed value of $\varepsilon'>0$ will be determined later),   there exists $   c >0$   such that, if     $ \disp (G) < c   $ then

  \begin{equation} \label{eqB}{ \dist   (T_x({\mathcal L}(x), F_x)) < \varepsilon'   } \ \ \forall \ x \in M.\end{equation}
\noindent Recall that $\delta_0, n_0$ (depending on $c$, which depends on $\varepsilon'$)  were defined by      Lemma~\ref{LemmaMane}. Therefore, each leaf ${f^n (\mathcal L(y) \bigcap B_{\delta_0}^n(x))}$ is part of a single leaf of a   foliation associated to a graph $  G_n$, for all $n \geq 0$. Besides, Lemma \ref{LemmaMane} states that
\begin{equation} \label{eqE}{ \disp    G_n < c \ \ \forall \, n \geq n_0 } \end{equation}
\noindent From   Inequalities (\ref{ocho}) and (\ref{eqE}) we deduce that
$$m ^{f^n (\mathcal L(y))} (f^n (\mathcal L(y) \cap B_{\delta_0}^n(x))) \leq m ^{f^n (\mathcal L(y))} (f^n (\mathcal L(y) )) \leq $$ $$[(1 + \disp G_n) \, \delta \,]^{\mbox{ \footnotesize dim }  {F}} < [(1 + c) \, \delta \,]^{\mbox{ \footnotesize  dim } {F}} \ \ \ \forall \ n \geq n_0. $$
Thus, there exists $K >0$ such that $$m ^{f^n (\mathcal L(y))} (f^n (\mathcal L(y) \cap B_{\delta_0}^n(x))) \leq K \ \ \ \forall \ n \geq 0.$$
So,  Assertion {\rm (c)}   of   Proposition \ref{PropositionMane} is proved for each fixed value of $\varepsilon'>0$.

Next, we prove {\rm (d)}.  From   Inequalities (\ref{eqB}) and (\ref{eqE}), we deduce that
  \begin{equation} \label{eqB'}\dist  \left(F_{f^n(x)}, T_{f^n(x)} f^n ({\mathcal L} (x))\right) < \varepsilon' \ \ \forall \ x \in M, \ \forall \ n \geq n_0. \end{equation}
   Finally, we fix  $\varepsilon'>0$ (depending on the  given value of $\varepsilon >0$), such that $0 < \varepsilon' < \varepsilon$ and such that for all $\dim F$-dimensional sub-bundles ${L}$ that satisfy $\dist  (L, F) < \varepsilon'$, the following inequality holds: \begin{equation}\label{eqn13ymedio} e^{-\varepsilon} \leq \frac{\left|\det df_{x}|_{{ L(x)}}\right|}{ \left|\det df_{x}|_{{ F_x}}\right| } \leq e^\varepsilon \ \ \forall \ x \in M.\end{equation}  Therefore, (\ref{eqB'}) implies:
   $$e^{-\varepsilon} \leq \frac{\left|\det df_{f^j (x)}|_{{ T_{f^j(x)} f^j{\mathcal L}(x)}}\right|}{\left|\det df_{f^j (x)}|_{{ F_{f^j(x)}}}\right|} \leq e^\varepsilon \ \ \forall \ x \in M, \ \ \forall \   j \geq n_0.$$
    The  latter inequality implies {\rm (d)}. Finally, {\rm (b)}   is obtained from (\ref{eqB'}) taking into account that $\varepsilon'$ was chosen smaller than $\varepsilon$. \hfill $\ \Box$

\vspace{.3cm}

\noindent\textsc{ End of the proof of Lemma \ref{Lemma 1}:}

For the given value of $\varepsilon >0$, we construct $\delta_0, K >0$ as in Proposition \ref{PropositionMane}. Consider any   finite partition $\alpha = \{A_h\}_{1 \leq h \leq k}$, where $k= \#(\alpha)$, such that $\diam  \alpha = \max_{h= 1\ldots k} \diam \{ A_h\} < \delta_0$.

 For each $A \in \alpha$  construct an open set $V_A \subset M$ also of diameter smaller than $\delta_0$, containing $A$. Construct a $\dim F$ local foliation ${\mathcal L}_A  $ in $V_A$ satisfying Proposition \ref{PropositionMane}.
Construct also a $C^1 $  submanifold $W_{A}$   transversal  to ${\mathcal L}_A$.
%\pause

Take $\alpha^n = \vee _{j= 0}^n f^{-j} (\alpha) = \{X_i\}_{1 \leq i \leq k_n}$, where $k_n = \#(\alpha^n)$. For all $   X_i \in \alpha^n $, there exists $A_{h_i} \in \alpha$ such that $\ X_i \subset A_{h_i}$. Denote ${\mathcal L}_i:= {\mathcal L}_{A_{h_i}}$ and $W_i := W_{A_{h_i}}$. %\pause
 Since ${\mathcal L}_i$ is $C^1$-trivializable, by Fubini's Theorem we have:

\begin{equation}\label{eqnFormula_m(C)} m(C) =   \sum_{i= 1}^{k_n} \int_{z \in W_{i} } d \mu^{W_{i}} \int _{y \in {\mathcal L}_i(z)}   {{\mathbf {1}}_{C \bigcap X_i} } \, \phi_{i}  \, dm^{   {\mathcal L}_i(z)} \ \ \ \ \ \forall \ C \in {\mathcal B},\end{equation}
where  ${\mathcal B}$ is the Borel sigma-algebra, ${\mathbf 1}_{C \bigcap X_i}$ is the characteristic function of the set $C \bigcap X_i$, and $\phi_{i}$ is a continuous function which depends on $A_{h_i} \in  \alpha$. Precisely, $\phi_i$ is the Jacobian of the $C^1$-trivialization of the foliation ${\mathcal L}_{i} $. So,  there are  at most $k = \#(\alpha)$ different local foliations ${\mathcal L}_i$, $k$ different continuous functions $\phi_i$, and $k$  different transversal manifolds $W_i$,   which allow Formula (\ref{eqnFormula_m(C)}) work  for any value of $n$ and for any $C \in {\mathcal B}$.
 %\pause

  Denote ${\widehat y} = f^n(y) \in f^n ({\mathcal L}_i(z) \bigcap X_i) =: {\mathcal L}_i^n (z)$:  \begin{equation}\label{equation3}\ m(C)=\sum_{i= 1}^{k_n} \int_{z \in W_{i} } d \mu^{W_{i}}  \int _{{\widehat y} \in {\mathcal L}_i^n(z) }  [{\mathbf {1}}_{C \cap X_i} \   \phi_{i} ](f^{-n}({\widehat y}))         \left|\det df^{-n}|_{T_{\widehat y}{\mathcal L}_i^n }    \right| \,   dm^{  {\mathcal L}_i^n(z)}.     \end{equation}

 By Part {\rm (d)}   of Proposition \ref{PropositionMane}: \begin{equation} \label{eqn4aa}\left|\det df^{-n}|_{T_{\widehat y}{\mathcal L}_i^n }  \right| \leq K e^{n \varepsilon} \, \left|\det df^{-n}|_{F_{\widehat y} }  \right|. \end{equation}
 Recall Formula (\ref{eqnPsiSubn}) defining $\psi_n(y)$.
 Since $f^n(y) = \widehat y$, we have
   $$\log\left|\det df^{-n}|_{F_{\widehat y} } \right|= \psi_n(f^{-n}(\widehat y)),$$
 which together with inequality (\ref{eqn4aa}) and equality (\ref{equation3}) gives:
 \begin{equation}  \label{equation4} \left|\det df^{-n}|_{T_{\widehat y}{\mathcal L}_i^n }  \right| \leq K e^{n \varepsilon}  e^{\psi_n (f^{-n}({\widehat y}))}\end{equation}
 \begin{equation}
 \label{equation3aa} m(C)\leq K e^{n \varepsilon} \sum_{i= 1}^{k_n} \int_{z \in W_{i} } d \mu^{W_{i}}  \int _{{\widehat y} \in {\mathcal L}_i^n(z) }  [{\mathbf {1}}_{C \cap X_i} \   \phi_{i} ](f^{-n}({\widehat y}))       \,  e^{\psi_n (f^{-n}({\widehat y}))} \,   dm^{  {\mathcal L}_i^n(z)}.
 \end{equation}
  \ By Riesz Representation  Theorem there exists a finite measure $\nu_n$ such that

 $$\int \!{h} \; d \nu_n \!=  \! \sum_{i= 1}^{k_n} \!\int_{z \in W_{i} }\!\! d \mu^{W_{i}}  \!\!\int _{{\widehat y} \in {\mathcal L}_i^n(z) }  \!\!\!\! {[({\mathbf{1}}_{X_i} \cdot \phi_i \cdot   h  ) \circ  f^{-n}] \, \, ({\widehat y})}  \;     dm^{    {\mathcal L}_i^n(z)}   \ \;  \forall \, h \in C^0(M, \mathbb{R}).  $$

From inequality (\ref{equation3aa})   and the above definition of $\nu_n$, we conclude
$$m(C) \leq K e^{n \varepsilon} \int {\mathbf {1}}_{C} \,  e^{\psi_n} \, d \nu_n \leq K e^{n \varepsilon} \int_C e^{\psi_n} \, d \nu_n.$$
Statement (ii) of Lemma \ref{Lemma 1} is proved.

Let us prove Statement (i). We must   show that there exists a constant $K_0 >0$, independent of $n$, such that $\nu_n(X) \leq K_0$ for all $X \in \alpha^n$, and for all $n \geq 0$.
In fact, recall that ${\mathcal L}_i^n(z) = f^n({\mathcal L}_i(z) \cap X_i) \subset f^n({\mathcal L}_i(z) \cap B_{\delta_0}(y))$ for all $z \in W_i$ and for all $ y \in  {\mathcal L}_i(z) \cap X_i$. Thus, applying Property  {\rm (c)}   of Proposition \ref{PropositionMane} we have $$m^{   {\mathcal L}_{i}^n  }({\mathcal L}_{i}^n(z)) \leq K_1,$$ for some constant $K_1 >0$ which is independent of $n$. From the construction of the measure $\nu_n$:
$$\nu_n(X_i) = \int_{z \in W_{i} } d \mu^{W_{i}}  \int _{{\widehat y} \in {\mathcal L}_i^n(z) }  {[(\mathbf{1}_{X_i} \cdot \phi_{i}) \circ  f^{-n}] \, ({\widehat y}) }  \,     dm^{    {\mathcal L}_i^n(z)} \leq $$ $$\mu^{W_{i}}(W_{i})  \,    \|\phi_i \|_{C^0} \, \, m^{   {\mathcal L}_{i}^n  }({\mathcal L}_{i}^n(z))  \leq K_1 \, \mu^{W_{i}}(W_{i})  \,    \|\phi_i \|_{C^0}.$$
Since the number of different local foliations ${\mathcal L}_i$   is equal to the number $k$ of pieces  of the given partition $\alpha$, which is independent of $n$, we obtain:
$$\nu_n(X_i) \leq K_1 \, \max_{A \in \alpha} \{\mu^{W_{A}}(W_{A})  \,    \|\phi_A \|_{C^0}\}  =: K_0  ,$$ where $K_0$ depends only on the partition $\alpha$ and not on $n$. \hfill $\ \Box$

\section{Proof of Lemma \ref{Lemma 2}} \label{sectionProofLemma2}

Choose $\{\varphi_i\}_{i \geq 1}$   dense in $C^0(M, [0,1])$ and define $\dist  ^*$ in ${\mathcal P}$:

\begin{equation}\label{3}\dist  ^* (\mu_1, \mu_2) := \left| \int \psi\,d \mu_1  - \int   \psi\, d \mu_2\right| + \sum_{i= 1}^\infty \frac{\left|\int \varphi_i\, d\mu_1 -  \int \varphi_i\, d \mu_2\right|}{2^i}.\end{equation}

By hypothesis,  a measure $\mu \in {\cal P}_f$ and two small numbers $\varepsilon >0$ and $\delta >0$   are arbitrarily given. We must construct an adequate finite partition $\alpha$ of $M$, with diameter smaller than $\delta$,   satisfying   Lemma \ref{Lemma 2}.

 Take $\delta_1 >0$ such that  $\dist  (x,y) < \delta_1$ $\Rightarrow$ $|\psi(x)- \psi(y)| < {\varepsilon}/{5}$.

  \noindent Take $\alpha$ such that $\diam (\alpha) \leq \min (\delta, \delta_1)$,  $\mu (\partial X)= 0 \ \forall \ X \in \alpha$.
  This construction implies
   \begin{equation} \label{eqn66a} \lim_{n \rightarrow + \infty} \mu_n(X) = \mu(X)\end{equation} for all $ X \in \alpha ^q = \bigvee _{j= 0}^q f^{-j}(\alpha)$, for all $  \ q \in \mathbb{N}$ and for all $ \{\mu_n\}_n \subset {\mathcal P} $ such that $  \lim^*_n\mu_n = \mu$.
   Also  \begin{equation}
   \label{eqn64}
   |\psi_n(y) - \psi_n(x)| \leq \sum_{j= 0}^{n-1} |\psi (f^j(y))  - \psi(f^j(x))|  \leq \frac{n \varepsilon}5 \ \ \forall \ x,y \in X , \forall  X \in \alpha^n.\end{equation}
\noindent Recall that $h_{\mu}(\alpha) := \lim_{q \rightarrow + \infty} H(\alpha^q, \mu)/q$, where $$H (\alpha^q, \mu) := -\sum_{X \in \alpha^q} \mu(X) \log \mu(X).$$
Fix $q \in \mathbb{N}^+$ such that $ H(\alpha^q, \mu)/q < h_{\mu}(\alpha)   + {\varepsilon}/5$. From (\ref{eqn66a}):

  \begin{equation}
  \label{eqn62}
  \lim_{n \rightarrow + \infty} \frac{H(\alpha^q, \mu_n)}q = \frac{H(\alpha^q, \mu)}q < h_{\mu}(\alpha)   + \frac{\varepsilon}5 \end{equation} for any sequence $\mu_n \in {\mathcal P}$ such that $\lim^*_n\mu_n = \mu$.

   %Choose $\{\varphi_i\}_{i \geq 1}$   dense in $C^0(M, [0,1])$ and define $\dist  ^*$ in ${\mathcal P}$:
%
%\begin{equation}\label{3}\dist  ^* (\mu_1, \mu_2) := \left| \int \psi\,d \mu_1  - \int   \psi\, d \mu_2\right| + \sum_{i= 1}^\infty \frac{\left|\int \varphi_i\, d\mu_1 -  \int \varphi_i\, d \mu_2\right|}{2^i}.\end{equation}

Using (\ref{eqn62}),  fix $0 <\varepsilon_0^* < \varepsilon/5$ such that
\begin{equation}
\label{eqn63}
\sigma \in {\mathcal P}, \ \dist  ^*(\sigma, \mu) \leq  \varepsilon_0^* \Rightarrow |H(\alpha^q, \sigma) - H(\alpha^q, \mu)| \leq \frac{q \varepsilon}5.\end{equation}
Such a value of $\varepsilon_0^*$ exists; otherwise we could construct a sequence of probability measures $\mu_n$  converging to $\mu$ and such that $|H(\alpha^q,  \mu_n) - H (\alpha^q, \mu)| > q \, \varepsilon/5$ for all $n \in \mathbb{N}$. This inequality contradicts the equality at left in (\ref{eqn62}).

For any fixed $0 <\varepsilon^* < \varepsilon^*_0$  we denote $C_n = C_n(\varepsilon^*)$  defined by equality (\ref{equationCn}).
  Consider $$\alpha^n\bigvee\{C_n\}:= \{X_i \cap C_n: \   \ X_i \in \alpha^n, \ \ X_i \cap C_n \neq \emptyset\}.$$ Denote $ k_n := \# (\alpha \bigvee \{C_n\})$.
For each $C_n \cap X_i \in \alpha^n \bigvee \{C_n\}$,   choose one  point $x_i \in C_n \bigcap X_i$. Consider the integral $I(\psi_n, C_n, \nu_n)$ defined by equality (\ref{eqn*}), and apply inequality (\ref{eqn64}):
\begin{equation} \label{0}I_n:= I(\psi_n, C_n, \nu_n) = \int_{C_n} e^{\psi_n} \, d \nu_n = \sum_{i=1}^{k_n} \int_{y \in C_n \cap X_i} e^{\psi_n(y)} d \nu_n(y) \leq $$ $$\sum_{i= 1}^{k_n} e^{n \varepsilon/5}   e^{{\psi_n(x_i)}} \nu_n(C_n \cap X_i). \end{equation}
\noindent By hypothesis $\nu_n(X) \leq K$ for all $X \in \alpha^n$. So:
 \begin{equation} \label{1}I_n \leq K e^{n \varepsilon/5} \sum_{i= 1}^{k_n} e^{\psi_n(x_i)}\end{equation}

  \noindent Define $p_i := e^{\psi_n (x_i)}/L$, where $ L:= \sum_{i= 1}^{k_n} e^{\psi_n (x_i)}$. Note that  $\sum_{i= 1}^{k_n} p_i = 1$.
Then:
\begin{equation} \label{8}\log \sum_{i= 1}^{k_n} e^{\psi_n (x_i)} = \sum_{i= 1}^{k_n} {\psi_n (x_i)} p _i - \sum_{i= 1}^{k_n} p_i \log p_i.\end{equation}
Taking logarithm in (\ref{1}) and using (\ref{8}), we obtain:
$$\log I_n \leq {\log K} + \frac{n \varepsilon}{5} + \sum _{i=1}^{k_n} \psi_n(x_i) p_i - \sum_{i=1}^{k_n} p_i \log p_i.$$
From Equalities (\ref{eqnPsi}) and (\ref{eqnPsiSubn}):
$$\sum_{i= 1}^{k_n} \psi_n(x_i) p_i = \sum _{i= 1}^{k_n} \sum_{j= 0}^{n-1}   \int p_i\, \psi \   d \delta_{f^j(x_i)},$$
and thus:
\begin{equation} \label{4}\log I_n \leq {\log K}  + \frac{n \varepsilon}5 + \sum_{i= 1}^{k_n} \sum_{j= 0}^{n-1} \int p_i \, \psi \  d \delta_{f^j(x_i)} -\sum_{i= 1}^{k_n} p_i  \log p_i.  \end{equation}
  Let $\sigma_{n,x}$ be the empirical probability according to  Definition \ref{definitionEmpirical}. We construct  $\mu_n   \in {\mathcal P} $ by the following equality:
\begin{equation} \label{7}\mu_n := \frac{1}{n} \sum_{i= 1}^{k_n} \sum_{j= 0}^{n-1}   p_i \  \delta_{f^j(x_i)} = \sum_{i= 1}^{k_n} p_i \sigma_{n,x_i}.\end{equation}
Since
$ x_i \in C_n $ we have that $\dist^*(\sigma_{n, x_i}, \mu) < \varepsilon^*$ - see equality (\ref{equationCn}). Since the $\varepsilon^*$-balls defined with the metric $\dist^*$ by equality (\ref{3}) are convex, and $\mu_n$ is a convex combination of the measures $\sigma_{n, x_i}$, we deduce  $$ \dist  ^*(\sigma_{n, x_i}, \mu) \leq \varepsilon^* \Rightarrow \dist  ^*(\mu_n, \mu) \leq \varepsilon^*. $$
From the construction of $\dist^*$ by equality (\ref{3}), we obtain $\left|\int \psi\, d\mu_n  - \int \psi \, d\mu \right| \leq \varepsilon^*< \varepsilon/5$. Therefore:  $$  \int \psi \, d \mu_n   \leq \int \psi \, d \mu + \frac{\varepsilon}5,  $$
which togethter  with (\ref{4}) and (\ref{7}) implies:
$$\log I_n  \leq
\log K + \frac{2n \varepsilon }5 + n\int \psi\, d \mu - \sum_{i= 1}^{k_n} p_i \log p_i.$$
In Assertion \ref{assertionEntropy} of the appendix we prove the following statement:

\em There exists $n_0 \geq 0$ such that
$$- \sum_{i=1}^{k_n} p_i \log p_i \leq \frac{n \, \varepsilon}5 + \frac{n}{q}  H(\alpha^q, \mu_n) \ \ \forall n \geq n_0$$ \em
 Therefore:
$$\log I_n  \leq
 \log K + \frac{3n \varepsilon }5 + n\int \psi\, d \mu + \frac{ nH(\alpha^q, \mu_n)}{q}.$$
 By the construction of  $\varepsilon_0^*$ in (\ref{eqn63}), and since $\dist  ^* (\mu_n, \mu) < \varepsilon^* < \varepsilon_0^* $, we deduce: $$ \left|\frac{H(\alpha^q, \mu_n)}q - \frac{H(\alpha^q, \mu)}q\right| \leq \frac{\varepsilon}5. $$
 So $$\log I_n  \leq
\log K +  \frac{4n \varepsilon }{5}  + n \, \int \psi\, d \mu + \frac{n\, H(\alpha^q, \mu)}q. $$
Finally, using the choice of $q$ by inequality (\ref{eqn62})  we conclude
   $$\log I_n  \leq
\log K +  n\, \varepsilon   + n \, \int \psi\, d \mu + n \, h_{\mu}(\alpha), $$
 ending the proof of Lemma \ref{Lemma 2}. \hfill $\Box$

\section{Appendix} \label{sectionAppendix}

In this appendix we check some technical assertions that were used in the proofs of Sections \ref{sectionProofLemma1} and \ref{sectionProofLemma2}.

%\vspace{.2cm}

\begin{assertion}
\label{assertionG*1}

Let  $\delta >0$ be   such that for all $x \in M$ $$\exp_x \colon \{v \in T_xM: \|v\| \leq 3 \delta\} \to B_{3 \delta}(x) \subset M $$ is a diffeomorphism. Let $n_0 >0$.

Then, there exists $0 <\delta_0< \delta$  such that
for all $ x \in M$, for all $ 0 \leq n \leq n_0$ and for any graph  $G$  \em (defined in ${\mathbf B}^{T_xM}_{\delta}({\mathbf 0}) \subset T_xM$)  \em with $$\mbox{\em disp } G < 1/2,$$
  there exists a  graph $G_n$ \em (defined in ${\mathbf B}^{T_{f^n(x)} M}_{\delta}({\mathbf 0})  $) \em  satisfying the following condition:

  \begin{equation}
  \label{eqn88} \begin{array}{l}
                  \mbox{ for all }   y = \exp_{x} (v_1 + v_2 + G(v_1, v_2)) \in B_{\delta_0}(x) \\
                  \mbox{ there  exists } (u_1, u_2) \in {\mathbf B}_{\delta}^{E_{f^n(x)}}({\mathbf 0})\times {\mathbf B}_{\delta}^{ F_{f^n(x)}}({\mathbf 0}) \\
                  \mbox{ where } u_1 \mbox{ depends only on } v_1 \\
                  \mbox{ and } f^n(y) = \exp_{f^n(x)}(u_1 + u_2 + G_n(u_1, u_2)).
                \end{array}
  \end{equation}
      %\begin{equation}  \mbox{ for all }   y = \exp_{x} (v_1 + v_2 + G(v_1, v_2)) \in B_{\delta_0}(x) \mbox{ there  exists } (u_1, u_2) \in {\mathbf B}_{\delta}^{E_{f^n(x)}}({\mathbf 0})\times {\mathbf B}_{\delta}^{ F_{f^n(x)}}({\mathbf 0})\end{equation}
%      \begin{equation} \label{eqn89}
%      \mbox{ where } u_1 \mbox{ depends only on } v_1, \mbox{ and }
%      \end{equation}
%      \begin{equation} \label{eqn90}
%  f^n(y) = \exp_{f^n(x)}(u_1 + u_2 + G_n(u_1, u_2)).\end{equation}
\end{assertion}

{\begin{figure}
\vspace{-.3cm}
\begin{center}\includegraphics[scale=.5]{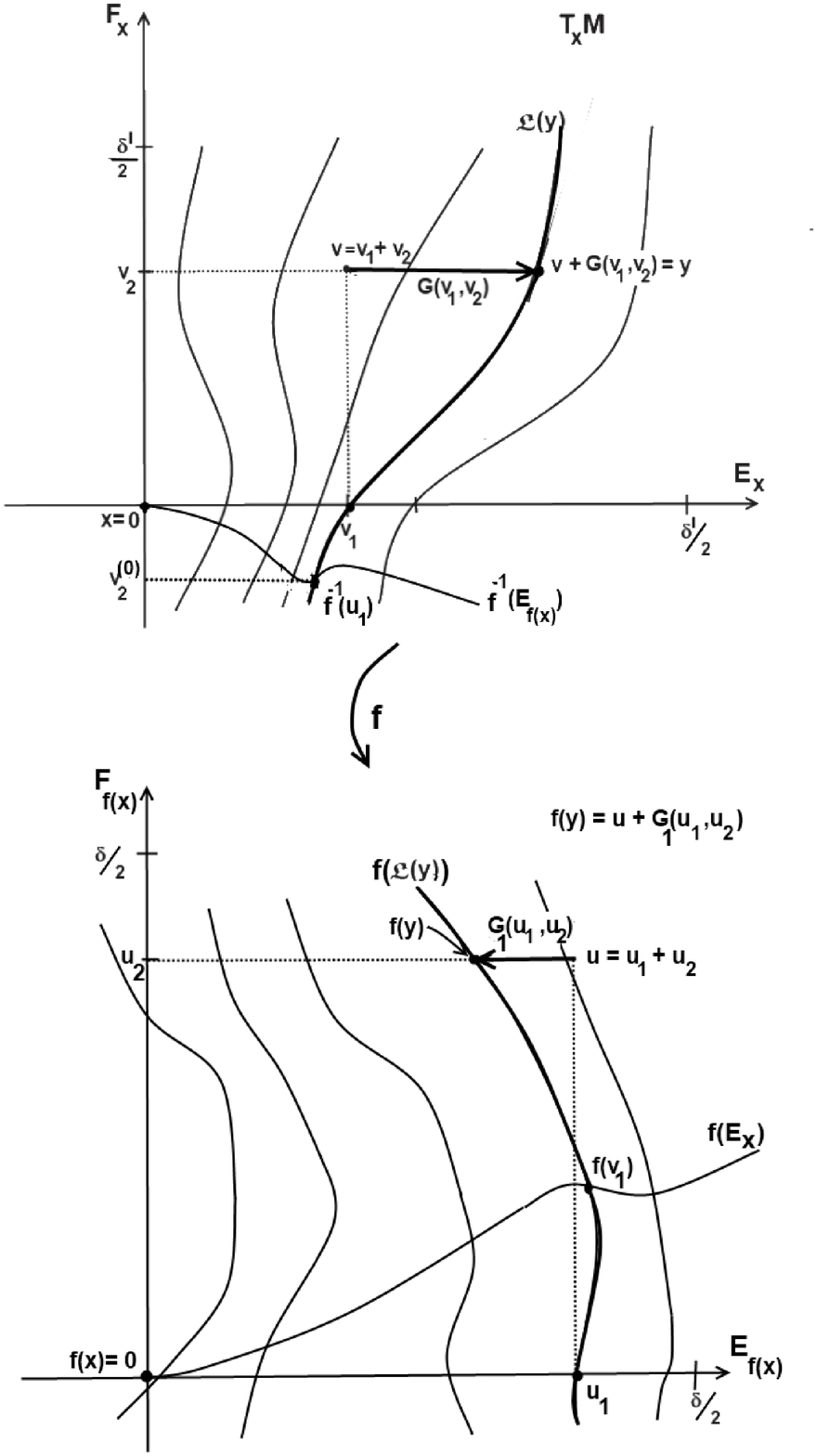}
\vspace{-1.6cm}
\caption{\label{Figure2} { \small The image $f({\mathcal L}(y))$ of the leaves ${\mathcal L}(y)$ near $x$ associated to the graph $G$, are associated to the graph $G_1$. (We omit the exponential maps $\exp_x$, $\exp_{f(x)}$).}}
\end{center}
\vspace{-1cm}
\end{figure}}

%\vspace{-1cm}

{\em Proof:} We will argue  by induction on $n \in \mathbb{N}$, to show that  for each $n \geq 1$, there exists $\delta_n >0$ and $G_n$ satisfying statement  (\ref{eqn88}). To prove Assertion \ref{assertionG*1}  it is enough to   take $\delta_0 :=\min\{\delta_1, \ldots, \delta_{n_0} \}$.

To simplify the notation along the proof, we will not write the exponential maps. From Definition \ref{definitiongraph}, recall the construction of the diffeomorphism $\Phi$ obtained from the graph $G$, which is a trivialization of the associated local foliation ${\mathcal L}$ (see the upper frame of Figure \ref{Figure2}).   Precisely, each leaf ${\mathcal L}(v_1)$ is obtained for constant $v_1 \in {\mathbf B}_{\delta}^{E_x}({\mathbf 0})$, and parametrized by $v_2 \in  {\mathbf B}_{\delta}^{F_x}({\mathbf 0})$   through  the formula $${\mathcal L}(v_1) : v_2 \mapsto  \Phi(v_1, v_2) := v_1 + v_2 + G(v_1, v_2), \ \mbox{ where } G(v_1, v_2) \in {\mathbf B}_{\delta}^{E_x}({\mathbf 0}).$$
Since $G(v_1, 0 )= 0$ we have $v_1 = \Phi(v_1, 0)$. i.e.  $$v_1 \in {\mathcal L}(v_1) \cap {\mathbf B}_{\delta}^{E_x}({\mathbf 0}).$$  Moreover, $$\pi_{F_x}(v_1 + v_2 + G(v_1, v_2)) = v_2 \mbox{ for all } v_2 \in {\mathbf B}_{\delta}^{F_x}({\mathbf 0}).$$ So
$$\pi_{F_x} {\mathcal L}(v_1) = {\mathbf B}_{\delta}^{F_x}({\mathbf 0}).$$

Besides, ${\mathcal L}(v_1)$ is uniformly transversal to $E_x$, for all $G$ with $\disp(G) < 1/2$. In fact  $$T_{v_1}{\mathcal L}(v_1) = \mbox{Im}(\mbox{Id}|_{F_x} + \partial G/\partial v_2),$$ the subspace $F_x$ is   transversal to $E_x$, and  $$\|\partial G/\partial v_2\| \leq \disp G < 1/2 .$$  Thus, since the leaf ${\mathcal L}(x)$ intersects $E_x$ at $\textsc{ 0}$, we deduce that there exists $0 <\delta' \leq \delta$, which is uniform for any $G$ with $\disp(G) < 1/2$, such that, if $\dist (x,y) < \delta '$ then $y'$ belongs to some leaf of the foliation ${\mathcal L}$. In other words, $$B_{\delta'}(x) \subset \mbox{Im}(\Phi),$$ i.e. there exists $(v_1, v_2) \in E_x \times F_x$ such that
$$\|v_1\| < \delta, \ \ \|v_2\| < \delta \mbox{ and } $$ $$ y= \Phi(v_1, v_2) = v_1 + v_2 + G(v_1, v_2) \ \ \mbox{ if } y \in B_{\delta'}(x).$$

Recall that, by Definition \ref{definitiongraph}, $\Phi$ is a diffeomorphism onto its image. Thus, for all $y \in B_{\delta'}(x)$, the point $\Phi^{-1}(y) = (v_1, v_2) \in {\mathbf B}_{\delta}^{E_x}({\mathbf 0})\times {\mathbf B}_{\delta}^{F_x}({\mathbf 0})$   depends $C^1$ on $y$.
We take $0 < \delta_1 < \delta' $ such that if $y  \in B_{\delta_1}(x)$, then
 \begin{equation} \label{eqn29}        \|\pi_{E_{f(x)}}\, f(y)\| < \delta/2, \ \ \| \pi_{F_{f(x)}} \, f    (y) \| < {\delta/2}.\end{equation}  Such a  value of $\delta_1 >0$ exists, and is independent of the graph $G$, because   $f$, \  $\pi_{E_x}$  and $\pi_{F_x}$    are uniformly continuous.

Taking  if necessary  a smaller value of $  \delta_1 $, the following two properties (A) and (B) hold for any graph $G$ with $\disp(G) < 1/2$ and for any $y   \in B_{\delta_1}(x) $:

  \textsc{ (A)}   The leaf $f({\mathcal L}(y)) $ intersects  ${ \mathbf{B}}_{\delta/2}^{E_{f(x)}} ( \mathbf{0}) \subset E_{f(x)}$ in a point  $u_1$  (see   Figure \ref{Figure2}).
  In other words \begin{equation}\label{eqn29bb}\mbox{ there exists } \ \ u_1 \in E_{f(x)}, \ \ \ \|u_1\| < \delta/2, \ \ \ f^{-1}(u_1) \in {\mathcal L}(y).\end{equation}

  %\vspace{-.9cm}

  \textsc{ (B)} The application
   $v_1 \in E_x \ \mapsto \ u_1 \in E_{f(x)}$   defined by (A)  for all $y = \Phi(v_1, v_2) \in B_{\delta_0}(x)$, is independent of $v_2$, and is a diffeomorphism onto its image.

  Property (A) is achieved due to the Implicit Function Theorem, since $f$ is a diffeomorphism, \ \ $\im df_x|_{E_x} = E_{f(x)}$,  and the local foliation ${\mathcal L}$ is  uniformly transversal to ${ \mathbf{B}}_{\delta}^{E_{x}} ( \mathbf{0}) \subset E_{x}$, while its leaf ${\mathcal L}(x)$    intersects $E_x$ at $0$.
 Property (B) is obtained because $f$ is a diffeomorphism  and the mapping $f(v_1) \in f(E_x) \mapsto u_1 \in E_{f(x)}$ is the holonomy along the leaves of the foliation $f({\mathcal L})$, which is $C^1$ trivializable and uniformly transversal to both $f(E_x)$ and $E_{f(x)}$. (See   Figure \ref{Figure2}.)

\vspace{.3cm}

 Let us show that the graph $G_1$ exists in $T_{f(x)} M$  satisfying Definition \ref{definitiongraph}   and Assertion  (\ref{eqn88})  for all $y \in B_{\delta_1}(x)$. We write $y = \Phi (v_1, v_2)  = v_1 + v_2 + G(v_1, v_2).$
   We have already determined $u_1 \in E_{f(x)}$ as a  diffeomorphic function of $v_1$, which    does not depend on  $v_2$. Let us determine $  u_2  \in F_{f(x)}$ and $ G_1(u_1, u_2) \in E_{f(x)}$  such that $f(y) = u_1 + u_2 + G_1(u_1, u_2)$ according to Figure \ref{Figure2}.  Consider the equation:  \begin{equation} \label{eqna3}
   f (v_1  + v_2  + G (v_1, v_2 )) = u_1 + u_2 + G_1(u_1 , u_2), \ \ \end{equation} $$ \mbox{ where }  \  u_2  \in F_{f(x)} \ \mbox{ and } \   G_1(u_1, u_2) \in E_{f(x)}.    $$
 Equation (\ref{eqna3}) is solved by \begin{equation} \label{eqn91} u_2 := \pi_{F_{f(x)}} f (v_1 + v_2 + G(v_1, v_2))  \in F_{f(x)},\end{equation} \begin{equation} \label{eqn92} G_1(u_1, u_2) := - u_1 + \pi_{E_{f(x)}} f (v_1 + v_2 + G(v_1, v_2))  \in E_{f(x)}.\end{equation}  The application $\Psi$ defined by $\Psi (v_1, v_2) = (u_1, u_2)$, where $u_1 $ and $ u_2$ are constructed as above,   is of $C^1$ class. In fact, $u_1$ depends only on $v_1$,   the mapping $v_1 \mapsto u_1$ is a diffeomorphism onto its image, and $u_2$ is constructed by Formula (\ref{eqn91}).

 Moreover,    $\Psi$  is $C^1$ invertible. In fact, on the one hand, we have that the application $v_1 \mapsto u_1$ is $C^1$ invertible and independent of $v_2$. On the other hand, for constant $v_1$ let us show that  Formula (\ref{eqn91}) applies $v_2 \mapsto u_2$  \   $C^1$-diffeomorphically. Precisely,  $G(v_1, 0) = 0$, \ $v_2 \in F_x$ \ $u_2 \in F_{f(x)}$, $G(v_1, v_2) \in E_{x}$ and $df_x\colon F_x \to F_{f(x)}$ is invertible.
 Thus, $\pi|_{F_{f(x)}} \, d f = df|_{F_x} \, \pi_{F_x}$. Taking derivatives in   equality (\ref{eqn91}) with respect to $v_2$ with constant $v_1$, we obtain:
 $$ \frac{\partial u_2 }{\partial v_2}   = df|_{F_x} \cdot \pi_{F_x} \left (Id|_{F_x} + \frac{\partial G(v_1, v_2)}{\partial v_2} \right ) = df|_{F_x} = ( df^{-1}|_{F_{f(x)}} )^{-1}. $$
 The second equality holds because   $G(v_1, v_2) \in E_x$ for all $(v_1, v_2)$, and so, the projection by $\pi_{F_x}$ composed with   any derivative of $G$, equals zero.
 We have proved that   $ {\partial u_2 }/{\partial v_2}$
 is invertible, and besides
  \begin{equation}\label{eqnfea}
   \left(\frac{\partial u_2}{\partial v_2} \right)^{-1} = \frac{\partial v_2}{\partial u_2} = df^{-1}|_{F_{f(x)}}
  \end{equation}
   concluding that the application  $\Psi$   is a diffeomorphism onto its image.

   Now, we define the mapping $\Phi_1$   by $$\Phi_1 (u_1, u_2) = u_1 + u_2 + G_1(u_1, u_2).$$   $\Phi_1$ is a $C^1$ diffeomorphism onto its   image, because its inverse is $ \Psi \circ \Phi \circ   f^{-1}$. So $G_1$ is $C^1$, and $\Phi_1$ is the $C^1$ trivialization of its associated foliation, which is, by construction, $f({\mathcal L}).$

     Finally,     (\ref{eqn29}) , (\ref{eqn29bb})  and (\ref{eqn92}) imply   $$\|G_1\| \leq \|u_1\| + \|\pi_{E_{f(x)}}f(y)\| < \delta/2 + \delta/2 = \delta  \ \ \ \mbox{ and }$$    $$\Phi_1^{-1} f(B_{\delta_1}(x)) \subset {\mathbf{B}}^{E_{f(x)}}_{\delta/2} ({\mathbf 0}) \times { \mathbf{B}}^{F_{f(x)}}_{\delta/2} ({\mathbf 0}).$$  Thus, $G_1\colon \Phi_1^{-1}(f (  B_{\delta_1}(y))) \to { \mathbf{B}}^{E_{f(x)}}_{\delta} ({ \mathbf{0}})$ can be $C^1$ extended to  be a graph $G_1\colon { \mathbf{B}}^{E_{f(x)}}_{\delta} ({\mathbf 0}) \times { \mathbf{B}}^{F_{f(x)}}_{\delta} ({\mathbf 0}) \to { \mathbf{B}}^{E_{f(x)}}_{\delta} ({ \mathbf{0}}). $

We have completed the first step of the inductive proof, since we have proved the existence of $\delta_1 >0$ and of the graph $G_1$ satisfying  (\ref{eqn88}). Naturally, $\disp  \ G_1$ is not necessarily upper bounded by $1/2$. So, we can not   exactly repeat  the same argument to prove the inductive step.   We will instead prove that there exists a uniform constant $c_1 >0$ such that   \begin{equation}\label{eqn30} \mbox{if } \disp  G < 1/2 \mbox{ then } \disp G_1 < c_1. \end{equation} If we prove  inequality (\ref{eqn30}) for some constant $c_1$, then we can end the inductive proof as follows.

Assume that   for some $n \geq 0$ there are  $\delta_n, c_n >0$ and a graph $G_n$ defined in ${\mathbf B}^{T_{f^n(x)} M}_{\delta} ({\mathbf 0})  $ satisfying (\ref{eqn88})   for all $y \in B_{\delta_n}(x)$, and such that $\disp  G_n \leq c_n $ for any graph $  G $ with $\disp  G < 1/2$. Thus, we can repeat the   above proof, putting $\min(\delta_n, \delta), c_n$ and $ G_n$  in the roles of $\delta, 1/2 $ and $G$  respectively. We deduce that there are  $\delta_{n+1}, c_{n+1}>0$ and a graph $G_{n+1} := (G_n)_1$, defined in ${\mathbf B}^{T_{f^{n+1}(x)} M}_{\delta} ({\mathbf 0})  $, which satisfies  (\ref{eqn88})   for all $y \in B_{\delta_{n+1}}(x)$, and such that $\disp  G_{n+1} < c_{n+1}$ for any graph $G$ for which  $\disp  G_n < c_n$. Thus, $G_{n+1}$ satisfies (\ref{eqn88})  for all $G$ such that $\disp  G < 1/2$.  Therefore, the inductive proof will be completed  once we show inequality (\ref{eqn30}).

So, let us find a constant $c_1$ satisfying inequality (\ref{eqn30}). To find $c_1$  we  will bound  from above the term $\|\partial G_1(u_1, u_2) /\partial u_2\|$.
From (\ref{eqn92}), and taking into account that $$\pi_{E_{f(x)}} \cdot df|_x =    df|_{E_x} \cdot \pi|_{E_x},$$ we obtain:
$$  \frac{\partial G_1(u_1, u_2) }{\partial u_2}     =      df|_{E_x} \cdot \pi|_{E_x} \cdot \left (Id|_{F_x} + \frac{\partial G(v_1, v_2)}{\partial v_2}\right ) \cdot \frac{\partial v_2}{\partial u_2} =df|_{E_x} \cdot   \frac{\partial G(v_1, v_2)}{\partial v_2} \cdot \frac{\partial v_2}{\partial u_2}. $$
Applying (\ref{eqnfea}) and the definition of dispersion, we deduce
\begin{equation} \label{eqnPrincipal} \disp (G_1) \leq    \|df|_{E_x}\| \cdot \disp(G) \cdot \|df^{-1}|_{F_{f(x)}}\| .\end{equation}
Thus, inequality (\ref{eqn30}) follows taking
$$c_1 := \max_{x \in M} \{ \|df|_{E_x}\| \cdot   \|df^{-1}|_{F_{f(x)}}\|\},$$
ending the proof of  Assertion \ref{assertionG*1}. \hfill $\Box$

\vspace{.2cm}

\begin{assertion}
\label{assertionG*2} Let  $\delta >0$ be   such that for all $x \in M$ $$\exp_x \colon \{v \in T_xM: \|v\| \leq 3 \delta\} \to B_{3 \delta}(x) \subset M $$ is a diffeomorphism.  For all $0 <\delta_0< \delta $ there exists $0 < c' < 1/2$  satisfying the following property:

\vspace{.3cm}

Assume that $G$   is a Hadamard graph   defined in ${\mathbf B}^{T_{ x} M}_{\delta} ({\mathbf 0}) $ such that \em $$\disp  G < c'.$$\em
Assume that there exists    $n \in \mathbb{N}$  and a graph $G_n$   in ${\mathbf B}^{T_{f^n(x)}M}_{\delta} ({\mathbf 0}) \subset T_{f^n(x)} M$ such that
\begin{equation} \label{eqn36}
 \begin{array}{l}
   \mbox{for all }   y = \exp_x(v_1 + v_2 + G(v_1, v_2) )\in  B_{\delta_{ 0}}^{n}(x) \\
  \mbox{there exists }      (u_1, u_2 ) \in {\mathbf B}_{\delta}^{E_{f^n(x)}}({\mathbf 0}) \times {\mathbf B}_{\delta}^{ F_{f^n(x)}}({\mathbf 0}), \\
   \mbox{where } u_1 \mbox{ depends only on } v_1  \mbox{ and } \\
   f^n(y)= \exp_{f^n(x)}(  u_1 \! \!+  \!u_2 \! + \!  G_{n} \! (u_1, u_2)).
 \end{array}
\end{equation}
Then, the following inequality holds for all  $y = \exp_x(v_1 + v_2 + G(v_1, v_2) )\in  B_{\delta_{ 0}}^{n}(x)$:

    \em
   \begin{equation}
     \label{eqn37}
     \Big\| \displaystyle{\frac{\partial  G _{n}}{\partial u_2}(u_1, u_2)} \Big\|     \leq   \big\|df^n|_{E_x} \big\| \cdot    \disp  G    \cdot \big\| df^{-n} |_{F_{f^n(x)}}\big\|. \end{equation}

\end{assertion}

\noindent{\em Proof:} To simplify the notation, we do not write the exponential maps.

Equality (\ref{eqn36}) can be written as follows:
\begin{equation}
\label{eqn38}
 f^n (v_1 + v_2 + G (v_1, v_2)) =    u_1 \! \!+  \!u_2 \! + \!  G_{n}    (u_1, u_2), \end{equation}
  where   $$(v_1, v_2) \in E_x\times F_x,\ \ \ \    (u_1, u_2) \in  E_{f^n(x)}\times F_{f^n(x)},$$
  $$G \in E_{x}, \ \  \ \  G_n \in E_{f^n(x)} \ \  \mbox{ and }
  $$ $$y = v_1 + v_2 + G(v_1, v_2) \in B_{\delta_0}^n(x).$$ Then:
\begin{equation}
\label{eqn38bbb}u_2 = \pi_{ {F_{f^n(x)}}} \, f^n (v_1 + v_2 + G(v_1, v_2)),\end{equation}
\begin{equation}
\label{eqn38ccc}G_n (u_1, u_2) = - u_1 + \pi_{E_{f^n(x)}} \, f^n (v_1 + v_2 + G(v_1, v_2)).\end{equation}

Taking derivatives in equality (\ref{eqn38bbb}) with respect to $v_2$, with constant $v_1$, and noting that $\pi|_{F_{f^n(x)}} \cdot df = df|_{F_x} \cdot \pi_{F_x}$, we obtain:
$$\frac{\partial u_2}{\partial v_2} = df^n|_{F_x} \pi_{F_x}  (Id|_{F_x} + (\partial G/\partial v_2)) = df^n|_{F_x} = (df^{-n}|_{F_{f^n(x)}})^{-1}. $$
In the second equality above, we used that $G(v_1, v_2)  \in E_x$ for all $(v_1, v_2)$ (recall Definition \ref{definitiongraph} of Hadamard graphs). Since $df^n|_{F_x}$ is invertible,   the linear transformation $ {\partial u_2}/{\partial v_2}$ is also invertible, and

$$\left ( \frac{\partial u_2} {\partial v_2}   \right )^{-1}=     \frac{\partial v_2} {\partial u_2}     = df^{-n}|_{F_{f^n(x)}}. $$

Now, we take derivatives in equality (\ref{eqn38ccc}) with respect to $v_2$ with constant $v_1$. We recall  that, by hypothesis,  $u_1$ depends  only on $v_1$, but not on $v_2$. Besides, we note that $\pi|_{E_{f^n(x)}} \cdot df^n = df^n|_{E_x} \cdot \pi|_{E^x}$. Thus:
$$ \frac{\partial G_{n}}{\partial u_2} \cdot \frac {\partial u_2}{\partial v_2} =      df^n|_{E_x} \cdot     \frac{\partial G}{\partial v_2}.$$
Thus:
$$ \frac{\partial G_{n}}{\partial u_2}  =   df^n|_{E_x} \cdot \frac{\partial G}{\partial v_2}    \cdot \frac{\partial v_2}{\partial u_2} = df^n|_{E_x} \cdot \frac{\partial G}{\partial v_2} \cdot df^{-n} |_{F_{f^n(x)}}. $$
So, after Definition \ref{DefinitionDispersion} of $\disp(G)$ we deduce
\begin{equation} \nonumber
\left \| \frac{\partial G_{n}}{\partial u_2}  \right \| \leq   \left\|df^n|_{E_x} \right\| \cdot    \disp  G    \cdot \left\| df^{-n} |_{F_{f^n(x)}}\right\|,     \end{equation}
proving inequality (\ref{eqn37}). \hfill  $\Box$

 \begin{assertion}
\label{assertionEntropy}

    There exists $n_0 \geq 1$ such that
$$  -\sum_{i= 1}^{k_n} p_i \log p_i\leq \frac{n\varepsilon}5 + \frac{n H(\alpha^q, \mu_n)}q \ \ \ \ \ \forall \ n \geq n_0, $$
\em where $\varepsilon >0$, $\alpha$ is a finite partition,   $\alpha^n= \bigvee_{j= 0}^{n} f^{-j}(\alpha),$ with  $ k_n \leq \#(\alpha^n),$  and $$0 \leq p_i\leq 1,\ \ \ \ \ \ \ \sum_{i= 1}^{k_n}p_i = 1, \ \ \ \ \ \ \ \ \ \   \mu_n := \frac{1}{n} \sum_{j= 0}^{n-1} \sum_{i= 1}^{k_n} p_i \delta_{f^j(x_i)}$$ with $x_i \in X_i$, \ $X_i \in \alpha^n$  and $$H(\alpha^q, \mu_n) := -  \sum_{\displaystyle{A \in \alpha^q} } \mu_n(A) \log \mu_n(A) .$$
\end{assertion}

{\em Proof:}  Denote $k:= \# \alpha$. Construct the probability measure
  $\pi_n := \sum_{i=1}^{k_n} p_i \delta_{x_i}   $. Then  $\pi_n (X_i) = p_i \ \  \forall \ 1 \leq i \leq k_n$ and $$H(\alpha^n, \pi_n)= -\sum_{i= 1}^{k_n}  p_i \log p_i $$
Fix $0 \!\leq l  \! \leq q-1$. Since $\alpha^{n+l}$ is thinner than $\alpha^n$, we have $H(\alpha^n, \pi_n) \leq H(\alpha^{n+l}, \pi_n)$. Thus
\begin{equation}\label{eqn67a} -\sum_{i= 1}^{k_n}  p_i \log p_i \leq H(\alpha^{n+l}, \pi_n),  \end{equation}
 where  $$\alpha^{n+l} = \vee_{j=0}^{n+l} f^{-j}\alpha = \Big( \vee_{j= 0}^{l-1} f^{-j}\alpha\Big ) \vee \Big (f^{-l}\big(\vee_{j= 0}^{n} f^{-j}\alpha \big)  \Big).$$

 \noindent Besides, for any two partitions $\alpha$  and $\beta$, and for any probability measure $\nu$ we have $H(\alpha \vee \beta, \nu) \leq H (\alpha, \nu) + H (\beta, \nu)$. Therefore \begin{equation} \label{eqn67b} H(\alpha^{n+ l}, \!\pi_n) \!\leq
 \! \sum_{j= 0}^{l-1}\! H(\alpha, \!{f^*}^j \!\pi_n) + H(f^{-l}\! \alpha^n, \!\pi_n),\end{equation}
where  the operator $f^*\colon {\mathcal P} \to {\mathcal P}$ in the space of probability measures is defined by $f^* (\nu)(B) = \nu(f^{-1}(B))$ for any measurable set $B$.

\noindent Since $H(\alpha, \nu) \leq \log(\#(\alpha)) = \log k$ for any probability measure $\nu$, and since $0 \leq l < q$, from Inequalities (\ref{eqn67a}) and (\ref{eqn67b}), we obtain:
 $$-\sum_{i= 1}^{k_n}  p_i \log p_i \leq q \log k + H(\alpha^n, {f^*}^{l} \pi_n).$$
If $n \geq   ( 10 \, q \,     \log k )\, / \, \varepsilon$, then \begin{equation}
\label{yoquese}
-\sum_{i= 1}^{k_n}  p_i \log p_i \leq \frac{n \varepsilon}{10} + H(\alpha^n, {f^*}^{l} \pi_n).   \end{equation}
Now we write:  $n = Nq + s, \  0 \leq s \leq q-1$. We have
$$H(\alpha^n, {f^*}^{l} \pi_n) \leq \sum_{h= 0}^{N-1} H(\alpha^q, {f^*}^{hq + l} \pi_n) + \sum_{j= Nq} ^{Nq+ s} H(\alpha, {f^*}^{j+l}\pi_n)  $$  {Arguing as above: $$\sum_{j= Nq} ^{Nq+ s} H(\alpha, {f^*}^{j+l}\pi_n) \leq (s+1) \log k \leq q \log k\leq \frac{n \varepsilon}{10}.$$ So, inequality (\ref{yoquese}) implies:}
%Thus,
$$-\sum_{i= 1}^{k_n}  p_i \log p_i \leq \frac{n \varepsilon}5 + \sum_{h= 0}^{N-1} H(\alpha^q, {f^*}^{hq + l} \pi_n)   $$
Taking  all values of $l$ such that $0 \leq l \leq q-1$ and adding the above bounds, we deduce: \begin{equation}\label{10}-q\sum_{i= 1}^{k_n}  p_i \log p_i  \leq   \frac{n q\varepsilon}5 + \sum_{h= 0}^{N-1} \sum_{l= 0}^{q-1} H(\alpha^q, {f^*}^{hq + l} \pi_n)
 \leq \frac{nq \varepsilon}5 + \sum_{j= 0}^{n-1} H(\alpha^q, {f^*}^j \pi_n).\end{equation}  Recall that the entropy $H$ of a partition with respect to a convex combination of probabilities, is not smaller  than   the convex combination of the entropies with respect to each of the probabilities. Since $ \displaystyle{ \mu_n =  \frac1n  \sum_{j= 1}^{n-1} \sum_{i=1}^{k_n} p_i \delta_{f^j (x_i)} = \frac1n\sum_{j= 0}^{n-1} {f^*}^j  \pi_n,} $ we deduce     $\displaystyle{\frac1n \sum_{j= 0}^{n-1} H(\alpha^q, {f^*}^j \pi_n) \leq H(\alpha^q, \mu_n). }$
 Substituting in inequality (\ref{10}) we conclude
$$\displaystyle{-q\sum_{i= 1}^{k_n}  p_i \log p_i
 \leq \frac{nq \varepsilon}5 + n H(\alpha^q, \mu_n),}$$  ending the proof of Assertion \ref{assertionEntropy}. \hfill $ \Box$

\vspace{1cm}

\noindent\textsc{ Acknowledgements } The authors thank Rafael Potrie and the anonymous referees   for their useful suggestions  and comments.
E.C. and M.C. were partially financed by Comisi\'{o}n Sectorial de Investigaci\'{o}n Cient\'{\i}fica of Universidad de la Rep\'{u}blica.   E.C. was also partially financed by Agencia Nacional de Investigaci\'{o}n e Innovaci\'{o}n (Uruguay).

\end{document}